\documentclass[preprint]{elsarticle}

\usepackage{lineno,hyperref}

\usepackage{color,amsmath,latexsym,amssymb,enumerate,bbm,cancel}
\usepackage[ruled,linesnumbered]{algorithm2e}
\usepackage{subfigure}
\usepackage{soul}
\modulolinenumbers[5]

\journal{Journal of Computational and Applied Mathematics}

\newcommand{\vect}[1]{\ensuremath{\mathbf{#1}}}

\newtheorem{theorem}{Theorem}
\newtheorem{lemma}{Lemma}

\newcommand{\R}{\mathbb{R}}

\renewcommand{\(}{\left(}
\renewcommand{\)}{\right)}

\newcommand{\ones}{\mathbbm{1}}
\newcommand{\zeros}{\mathbf{0}}

\newcommand{\flam}[1]{\ensuremath{f^{\left\{#1\right\}}}}

\newcommand{\lam}[2]{\ensuremath{{\vect{#1}}^{\left\{#2\right\}}}}

\newcommand{\lamT}[2]{\ensuremath{{{\vect{#1}}^{\left\{#2\right\}\intercal}}}}

\newcommand{\slam}[2]{\ensuremath{{{#1}^{\left\{#2\right\}}}}}

\newcommand{\slamc}[3]{\ensuremath{{#1}^{\left\{{#2}\right\}}_{#3}}}

\begin{document}
\begin{frontmatter}

\title{Relaxed Multirate Infinitesimal Step Methods for Initial-Value Problems}

\author[lbnladd]{Jean M. Sexton}
\ead{jmsexton@lbl.gov}
\author[smuadd]{Daniel R. Reynolds\corref{cor1}}
\ead{reynolds@smu.edu}

\cortext[cor1]{Corresponding author}
\address[lbnladd]{Lawrence Berkeley National Laboratory,
  1 Cyclotron Rd, Berkeley, CA 94720}
\address[smuadd]{Department of Mathematics, Southern Methodist
  University, Dallas, TX 75275-0156}

\begin{abstract}
This work focuses on the construction of a new class of fourth-order
accurate methods for multirate time evolution of systems of ordinary
differential equations.  We base our work on the Recursive Flux
Splitting Multirate (RFSMR) \cite{Schlegel2012,Schlegel2012a,Schlegel2012b}
version of the Multirate Infinitesimal Step (MIS) methods
\cite{Knoth2012,Knoth2014,Knoth1998,Schlegel2009,Wensch2009},
and use recent theoretical developments for Generalized Additive
Runge-Kutta methods \cite{Gunther2013,Sandu2015} to propose our
higher-order \emph{Relaxed Multirate Infinitesimal Step} extensions.
The resulting framework supports a range of attractive properties for
multirate methods, including telescopic extensions, subcycling,
embeddings for temporal error estimation, and support for changes to the
fast/slow time-scale separation between steps, without requiring any
sacrifices in linear stability.  In addition to providing rigorous
theoretical developments for these new methods, we provide
numerical tests demonstrating convergence and efficiency on a suite of
multirate test problems.
\end{abstract}


\end{frontmatter}



\section{Introduction}
\label{ch:introduction}

Increasingly, computational science requires large-scale simulations
that consistently and accurately couple distinct physical processes.
While the mathematical models for individual processes often have a
well-known type (hyperbolic, parabolic, etc.) and are suitable for
classical numerical integrators, the same cannot be said for the
coupled models.  These multiphysics models are often of mixed type,
may have limited differentiability, and involve processes that evolve
at dissimilar rates.  As such, many multiphysics simulations require
more flexible time integrators that may be tuned for these complex
problems.

In this work, we consider systems of autonomous ordinary differential
equations that may be characterized by two distinct dynamical time
scales, one ``fast'' and the other ``slow'',
\begin{align}
\label{eq:additive_form}
  y'(t) &= \flam{s}(y)+\flam{f}(y), \quad t_0\le t\le t_f,\\
  \notag
  y(t_0) &= y_0 \in \R^n.
\end{align}
We note that this \emph{additive} form of the problem
does not prohibit multirate applications where the solution variables
are also partitioned between the two scales,
\begin{align}
  \notag
  &\begin{bmatrix} y_{s}'\\ y_{f}' \end{bmatrix} =
    \begin{bmatrix}
      \overline{f}_{s}\(y_{s},y_{f}\)\\
      \overline{f}_{f}\(y_{s},y_{f}\)
    \end{bmatrix}, \quad   t_0\le t\le t_f\\
  \notag
  &y_{s}(t_0) = y_{s,0}\in\R^m, \quad
   y_{f}(t_0) = y_{f,0}\in\R^{n-m}\\
  \label{eq:partitioned_reformulation}
  \Rightarrow &\\
  \notag
  &y = \begin{bmatrix} y_{s}\\ y_{f} \end{bmatrix},
  \quad
  \flam{s}(y) = \begin{bmatrix}
    \overline{f}_s\\
    0
  \end{bmatrix},
  \quad
  \flam{f}(y) = \begin{bmatrix}0\\
     \overline{f}_f
  \end{bmatrix}.
\end{align}
Multirate algorithms evolve the initial-value problem
\eqref{eq:additive_form} using different time step sizes for the fast
and slow processes.  For such methods to result in efficiency
improvements over more traditional single-rate methods that evolve all
processes at the fast time scale, the slow right-hand side function
$\flam{s}$ must be significantly more costly than $\flam{f}$.


In this work, we consider multirate methods where
the slow operator is integrated explicitly and the fast operator
is integrated either explicitly or implicitly.  To simplify our
presentation, we denote the characteristic time-scales of
the slow and fast components of the multirate problem as
$\overline{dt}_{s}$ and $\overline{dt}_{f}$, respectively.  We define
the \textit{time-scale separation} between these as
$\overline{m}=\overline{dt}_{s}/\overline{dt}_{f}$.
We integrate these components using the time-step sizes ${dt}_{s} = h$
and ${dt}_{f}=h/m$, where $m$ is the smallest integer satisfying
$\overline{m} < m$.  We note that other authors weaken this
restriction, by allowing for variable step sizes at the fast time
scale \cite{Bremicker-Trubelhorn2017}.  While the methods proposed
here may be easily extended to such situations, the focus of this
paper is on the construction of higher-order methods; investigations of
adaptive fast substepping are left for future work.

As with single-rate problems, the order of accuracy of a multirate
method describes how the error behaves as $h\to0$.  Within the
multirate context, this overall error arises from multiple sources:
the error in the evolution of the fast components, the error in the evolution
of the slow components, and the coupling error between the two scales.
For a variety of multirate methods, the fast component method is a
composition of $m$ steps of a fast base method, while the slow
component method is one step of a slow base method -- in this case the
multirate method order depends on the order of both base methods and
on the coupling order.

\subsection{Runge--Kutta Theory for Multirate Methods}
In recent work, Sandu and G{\" u}nther introduced a theory for
Generalized-Structure Additively-Partitioned Runge-Kutta (GARK) methods
\cite{Sandu2015,Sandu2013gark}.  As the name suggests, this theory
may be applied to understand a wide range of Runge-Kutta-like time
integration methods, and notably provides both order conditions and a
linear stability theory for methods in this general format.  For
simplicity, GARK methods are formulated for autonomous systems of ODEs
in additively-split form,
\[
   \frac{dy}{dt} = f(y) = \sum_{q=1}^N f^{\left\{ q\right\} }\(y\).
\]
Here, we consider the case of $N=2$, corresponding to our target
problems \eqref{eq:additive_form}, i.e.
\[
   \frac{dy}{dt} = f^{\{f\} }\(y\) +
                   f^{\{s\} }\(y\).
\]
Here, GARK methods are uniquely determined by coefficients in an
expanded Butcher tableau,
\begin{align}
\begin{array}{c|c}
\lam{A}{f,f} & \lam{A}{f,s}\\ \hline
\lam{A}{s,f} & \lam{A}{s,s}\\ \hline\hline
\lamT{b}{f} & \lamT{b}{s}
\end{array} \label{eq:GARK}
\end{align}
where internal stages and solution update are given by
\begin{align}
  \notag
  \slamc{k}{f}{j} &= y_{n} +
                     h \sum_{l=1}^{s^{\{f\}}} \slamc{a}{f,f}{j,l} f^{\{f\}}\(\slamc{k}{f}{l}\) +
                     h \sum_{l=1}^{s^{\{s\}}} \slamc{a}{f,s}{j,l} f^{\{s\}}\(\slamc{k}{s}{l}\), \\
  \label{eq:GARK_step}
  \slamc{k}{s}{i} &= y_{n} +
                     h \sum_{l=1}^{s^{\{f\}}} \slamc{a}{s,f}{i,l} f^{\{f\}}\(\slamc{k}{f}{l}\) +
                     h \sum_{l=1}^{s^{\{s\}}} \slamc{a}{s,s}{i,l} f^{\{s\}}\(\slamc{k}{s}{l}\), \\
  \notag
  y_{n+1} &= y_{n} + h \sum_{l=1}^{s^{\{f\}}} \slamc{b}{f}{l} f^{\{f\}}\(\slamc{k}{f}{l}\) +
                    h \sum_{l=1}^{s^{\{s\}}} \slamc{b}{s}{l} f^{\{s\}}\(\slamc{k}{s}{l}\),
\end{align}
in which the slow stage indices range over $i=1,\dots,s^{\{s\}}$ and the
fast stage indices range over $j=1,\dots,s^{\{f\}}$.  In particular,
we note that this constitutes a ``hybrid'' of partitioned Runge-Kutta
(PRK) and additive Runge-Kutta (ARK) formulations: there are distinct
stages for each right-hand side component (PRK-like), and all
right-hand side functions are used to update each stage (ARK-like).

As is typical for additive Runge-Kutta methods, Sandu and G{\"u}nther
make the simplifying assumption of \emph{internal consistency}, i.e.
\begin{equation}
  \label{eq:consistency_conditions}
  \slamc{c}{q}{i} = \sum_{j=1}^{s^{\{f\}}} \slamc{a}{q,f}{i,j} = \sum_{j=1}^{s^{\{s\}}} \slamc{a}{q,s}{i,j}\qquad
  \forall\, i=1,\ldots,s^{\{q\}}, \quad q=f,s.
\end{equation}
Order conditions for methods up to fourth order (and for GARK tableau
of arbitrary size) are provided in matrix-vector form in
\cite{Sandu2015,Sandu2013gark}, and additionally in elementwise form
in \cite{Sandu2013gark}.  As our work focuses on the derivation of a
new class of fourth order methods, we reproduce their matrix-vector form
here, where we assume internal consistency
\eqref{eq:consistency_conditions}:
\begin{align}
  \lamT{b}{\sigma}\ones^{\{\sigma\}} &= 1, \label{eq:mvOrderGARK1}\\
  \lamT{b}{\sigma}\lam{c}{\sigma} &= \frac12, \label{eq:mvOrderGARK2}\\
  \lamT{b}{\sigma}\(\lam{c}{\sigma}\times\lam{c}{\sigma}\) &= \frac13, \label{eq:mvOrderGARK3a}\\
  \lamT{b}{\sigma}\lam{A}{\sigma,\nu}\lam{c}{\nu} &= \frac16, \label{eq:mvOrderGARK3b}\\
  \lamT{b}{\sigma}\(\lam{c}{\sigma}\times\lam{c}{\sigma}\times\lam{c}{\sigma}\) &= \frac14, \label{eq:mvOrderGARK4a}\\
  \(\lam{b}{\sigma}\times\lam{c}{\sigma}\)^{\intercal}\lam{A}{\sigma,\nu}\lam{c}{\nu} &= \frac18, \label{eq:mvOrderGARK4b}\\
  \lamT{b}{\sigma}\lam{A}{\sigma,\nu}\(\lam{c}{\nu}\times\lam{c}{\nu}\) &= \frac{1}{12}, \label{eq:mvOrderGARK4c}\\
  \lamT{b}{\sigma}\lam{A}{\sigma,\mu}\lam{A}{\mu,\nu}\lam{c}{\nu} &= \frac{1}{24}. \label{eq:mvOrderGARK4d}
\end{align}
In these formulas, $\sigma$, $\nu$ and $\mu$ each range over the
values $\{f,s\}$, and we denote a column vector of ones in
$\R^{s^{\{\sigma\}}}$ as $\ones^{\{\sigma\}}$.  Throughout this
manuscript we denote standard matrix and vector multiplication with
adjacent objects (e.g., $\vect{b}^{\intercal}\vect{c}$
is an inner product), and we use the $\times$ operator or exponents to
denote component-wise multiplication (e.g., for $\vect{b},\vect{c}\in\R^n$ then
$(\vect{b}\times\vect{c})_i = b_ic_i$ and $(\vect{b}^2)_i = b_ib_i$
for $i=1,\ldots,n$).  We note that of the above equations,
\eqref{eq:mvOrderGARK1} corresponds to the order 1 conditions,
\eqref{eq:mvOrderGARK2} corresponds to the order 2 conditions,
\eqref{eq:mvOrderGARK3a}-\eqref{eq:mvOrderGARK3b} correspond to the order 3 conditions, and
\eqref{eq:mvOrderGARK4a}-\eqref{eq:mvOrderGARK4d} correspond to the order 4 conditions.
Hence for a two-component splitting, a fourth order method requires 28
order conditions to be met.  As with \cite{Gunther2016}, we
categorize these order conditions into two groups: all conditions that
include $\lam{b}{f}$ (i.e., the 14 conditions
\eqref{eq:mvOrderGARK1}-\eqref{eq:mvOrderGARK4d} with $\sigma=f$) are
considered as ``fast order conditions'', and the rest based on
$\lam{b}{s}$ are considered as ``slow order conditions.''

\subsection{Multirate Background}

Multirate methods and the understanding of their order and stability properties
have evolved over time; however the GARK framework can be used to
describe and analyze the full history of multirate methods having
Runge--Kutta type \cite{Gunther2016}.
In this section, we first give a summary of the
historical context of multirate methods, and in the following section
we discuss specifics of how to apply GARK theory to these methods.

The first multirate methods were introduced by Gear and Wells in the
early 1980's \cite{Gear1980, Wells1982, Gear1984}; these methods
integrated only a single scale at a time, using
linear interpolation to provide data from one time scale to the other.
This work was extended by Constantinescu and Sandu
to construct second order multirate methods based on explicit linear
multistep base methods \cite{Constantinescu2007,Sandu2009}.

A more general approach to construction of third-order multirate
methods was introduced in 1999 by Kv{\ae}rn{\o} \cite{Kvaerno1999}, and
was further elaborated by many authors
\cite{Bartel2002,Estep2012,Gunther2001,Kvaerno2000}.  In these
approaches, one assumes a partitioned formulation with variables
split into active (fast) and latent (slow).  From here, higher-order
interpolants between the two components are obtained by interleaving
the interpolation processes between the fast and slow scales.  A
particular strength of this work is their derivation of a set of
simplifying assumptions, so that once these are satisfied and the base
method is third order, the resulting multirate method will also be
third order.  More recently, Fok and Rosales constructed
\emph{linearly} fourth-order multirate methods using higher-order
interpolation between the slow and fast stages of a
particular Runge-Kutta base method \cite{Fok2016}.

An alternate approach for higher-order multirate schemes has been
pursued through extrapolation of lower-order methods.  The first such
work is from Engstler and Lubich's 1997 paper \cite{Engstler1997},
that generated higher-order solutions by extrapolating first-order
methods.  The efficiency of these methods was subsequently improved
through use of dense output formulas \cite{Engstler1997, Hairer1993,
  Hairer1990}. More recent methods of this type, including explicit
fast and explicit slow integration (explicit/explicit),
as well as implicit fast and explicit slow integration
(implicit/explicit), have been explored by Constantinescu and Sandu
\cite{Constantinescu2010,Constantinescu2013}, who investigated
extrapolation from first-order accurate methods.
These multirate extrapolation approaches increase the possible
order of the method at the cost of repeated evaluations of
the same time step.  Hence, although these methods are reasonable for
improvement one or two orders of accuracy, they may become
cost-prohibitive when generating higher-order methods.

Each of these approaches has included some effort to develop methods
which enforce desirable properties based on the originating PDE problem.
Several papers have been published which investigate the coupling of
finite element discretizations with multirate schemes \cite{Trahan2012,
Seny2013, Giraldo2008}. Constantinescu and Sandu introduce methods which
focus on preserving conservation properties for PDEs \cite{Constantinescu2007}.
Schlegel also posited some mass-preserving methods \cite{Schlegel2012}, that
focus on conservation of linear invariants in the multirate scheme.

In \cite{Gunther2016}, G{\"u}nther and Sandu not only analyzed many of
the aforementioned multirate methods using GARK theory, but also
proposed a general infrastructure for constructing multirate GARK
(MrGARK) methods.  In a recent preprint \cite{Sarshar2018}, Sarshar
and collaborators utilize this framework to develop a set of
fourth-order explicit/explicit multirate GARK methods, through
construction of the coefficient matrices $\lam{A}{f,s}$ and $\lam{A}{s,f}$ to
directly address the coupling conditions between two fourth-order base
methods.  That work additionally proposes a novel approach for
multirate time adaptivity, through construction of various sets of
embedding coefficients to separately measure the errors at the slow
and fast time scales, as well as the coupling error between scales.

Our work builds most closely off of recent work in Multirate
Infinitesimal Step (MIS) methods.  Originally developed by Wensch, Knoth and
Galant, MIS methods were constructed as a generalization of
split-explicit methods in the context of numerical weather prediction
\cite{Wensch2009,Klemp1978,Skamarock1992}.
A key contribution from \cite{Wensch2009} was the development of
a systematic approach to the order conditions for split-explicit methods
based on PRK theory \cite{Wensch2009}, allowing the development of
second and third order MIS methods for a variety of applications
\cite{Schlegel2012,Schlegel2012a,Schlegel2012b,Knoth2012,Knoth2014,Knoth1998,Schlegel2009,Wensch2009}.
The basic approach of these methods is that the slow right-hand side
function $\flam{s}$ is evaluated to provide forcing terms to create a
sequence of modified initial-value problems at the fast time scale,
\begin{align}
  \label{eq:mis_fast_ode}
  &v'(\tau) = f^{\{f\}}(\tau, v) + \sum_{j=1}^i \alpha_{i,j}
  f^{\{s\}}(t_n+c_j h, Y_j^{\{s\}}),\\
  \notag
  &v(t_n+c_ih) = Y_i^{\{s\}},
\end{align}
where $\tau\in[t_n+c_ih,t_n+c_{i+1}h]$, $i$ ranges over the
stages at the slow time scale, and where the coefficients
$\alpha_{i,j}\in\R$ are defined based on the slow Butcher tableau.  In
particular, we highlight the development of the Recursive
Flux-Splitting Multirate (RFSMR) methods in \cite{Schlegel2012} based
on this MIS theory. RFSMR methods are a specialization of MIS methods
which leverage a recursive formulation of the differential equations
at different scales in order to improve parallel performance.

Very recently, Sandu has posted a preprint that also extends MIS
methods to higher order \cite{Sandu2018}.  In his proposed \emph{Multirate
Infinitesimal GARK} (MRI-GARK) approach, Sandu changes the formulation
of the modified initial-value problems at the fast time scale
\eqref{eq:mis_fast_ode} to instead create a time-dependent forcing
function based off of the slow right-hand side evaluations,
\begin{equation}
  \label{eq:mri_gark_fast_ode}
  v'(\tau) = f^{\{f\}}(\tau, v) + \sum_{j=1}^{i+1} \gamma_{i,j}(\tau)
  f^{\{s\}}(t_n+c_j h, Y_j^{\{s\}}),
\end{equation}
where now the time-dependent coefficients $\gamma_{i,j}(\tau)$ must be defined appropriately.
Through this generalization of MIS methods, he has been able to create
fourth-order multirate methods, including those with a limited amount
of implicitness at the slow time scale.

Also very recently, Roberts and collaborators have posted a preprint
that explores nonlinear solver approaches and linear stability for
implicit/implicit multirate methods based on Sandu's MRI-GARK approach
\cite{Roberts2018}.

\subsection{GARK representation of MIS methods}\label{sec:MISGARK}

Prior to introducing our proposed methods, we first summarize the
analysis of MIS methods using the GARK formalism, first shown in
\cite[Theorem 4]{Gunther2013}.  A later, more detailed description can be found
in \cite{Sandu2018}.
In that analysis the problem is
considered to be in autonomous form in order to simplify the
exposition; however, we note that MIS approach may be easily extended
to non-autonomous problems \cite{Knoth2014}.

MIS methods are typically constructed using a pair of ``base''
Runge-Kutta methods: a variation of the ``outer'' method
$\left\{A^O,c^O,b^O\right\}$ is applied to the slow time scale, and a
variation of the ``inner'' method $\left\{A^I,c^I,b^I\right\}$ is
applied to the fast time scale.  The outer method is an explicit
Runge-Kutta method with $s^O$ stages, with the requirement that
$c^O_1=0$, $c^O_i \le c^O_j$ for $i<j$, and $c^O_{s^O}\le 1$.
The MIS scheme makes no assumption about the structure of the inner
method, aside from the requirement that it is a one-step method.  In
general, MIS methods may be run so that the inner method takes several
``fast'' time steps to evolve between the slow stages, e.g.~over the
interval $\left[t_n+c^O_{i-1}h,t_n+c^O_{i}h\right], i=2,\ldots,s^O$.  As in
\cite{Gunther2013}, both here and in our new theoretical work in
Section \ref{ch:RMIS}, we consider the case when only a single step of
the inner method is taken to evolve between each slow stage.  Since
any sequence of $m$ steps of a $q$-order, one-step method may be
equivalently written as a single step of a corresponding $q$-order
method, this in no way limits the applicability of the current theory.
We do note, however, that if the abcissae of the outer method are not
evenly spaced, i.e.~$c^O_i-c^O_{i-1} \ne c^O_j-c^O_{j-1}$ for some $i\ne j$,
then the Butcher tables for the inner one-step methods corresponding
to the outer stages $i$ and $j$ will not be identical.

As shown in \cite[Theorem 4]{Gunther2013}, the GARK tables
corresponding to the MIS scheme are given by $\lam{b}{s} = b^O$,
$\lam{c}{s} = c^O$, $\lam{A}{s,s} = A^O$,
\begin{align}
  \label{eq:MIS_Aff}
  \lam{A}{f,f} &= \begin{bmatrix}
    c_2^O A^I & \zeros & \cdots & \zeros\\
    c_2^O \ones^{\{s^I\}} b^{I\intercal} & \(c_3^O - c_2^O\)A^I & \cdots & \zeros\\
    \vdots &  & \ddots & \\
    c_2^O\ones^{\{s^I\}} b^{I\intercal} & \(c_3^O - c_2^O\)\ones^{\{s^I\}} b^{I\intercal} & \cdots & \(1-c_{s^O}^O\)A^I
  \end{bmatrix},\\
  \label{eq:MIS_Asf}
  \lam{A}{s,f} &= \begin{bmatrix}
    c_2^O \vect{g}_2 b^{I\intercal} & \cdots &
    \(c_{s^O}^O-c_{s^O-1}^O\)\vect{g}_{s^O} b^{I\intercal} &
    \zeros \end{bmatrix},\\
  \label{eq:MIS_Afs}
  \lam{A}{f,s} &= \begin{bmatrix}
    c^I \vect{e}_2^{\intercal} A^O\\
    \vdots\\
    \ones^{\{s^I\}}\vect{e}_{i}^{\intercal} A^O + c^I \(\vect{e}_{i+1}-\vect{e}_{i}\)^{\intercal} A^O\\
    \vdots\\
    \ones^{\{s^I\}} \vect{e}_{s^O}^{\intercal} A^O + c^I \(\vect{b}^{O\intercal}-\vect{e}_{s^O}^{\intercal} A^O\)
  \end{bmatrix},\\
  \label{eq:MIS_bf}
  \lamT{b}{f} &= \begin{bmatrix}
    c_2^O b^{I\intercal} & \(c_3^O-c_2^O\) b^{I\intercal} & \cdots &
    \(1-c_{s^O}^O\) b^{I\intercal} \end{bmatrix},
\end{align}
where the block matrices and vectors comprising these are defined
using standard matrix-vector notation,
\[
  \vect{g}_i\in\R^{s^O},\quad\text{with}\quad
  \left[\vect{g}_i\right]_j = \begin{cases}
    0,& j<i\\
    1,& j>=i\end{cases},
\]
and $\vect{e}_i$ is the $i$-th elementary basis vector.  The
corresponding abcissae for the multirate method are then $\lam{c}{s} =
c^O \in \R^{s^O}$ and
\begin{align}
  \label{eq:MIS_cf}
  \lam{c}{f} &= \begin{bmatrix} c_2^O c^I\\
    c_2^O\ones^{\{s^I\}} + \(c_3^O - c_2^O\) c^I\\
    \vdots\\
    c_{s^O}^O\ones^{\{s^I\}} + \(1-c_{s^O}^O\) c^I
  \end{bmatrix}.
\end{align}
We note that for MIS methods the number of slow stages matches the
outer table, $s^s = s^O$, and when each of these slow stages uses the
same inner table, $T_I$, the number of fast stages equals the
product of these stage numbers, $s^f = s^Os^I$.  We also note the
slight difference in presentation of $\lam{A}{f,s}$ from that
shown in \cite[Theorem 4]{Gunther2013}; in equation \eqref{eq:MIS_Afs}
the generic entry denoted by $i$ corresponds to the $i$-th block
row of the matrix.
Based on these Butcher tables for the GARK coefficients, G{\"u}nther
and Sandu prove a number of particularly beneficial properties of MIS
methods \cite{Gunther2013,Sandu2013multirate}:
\begin{itemize}
\item[(i)] The coefficients \eqref{eq:MIS_Aff}-\eqref{eq:MIS_cf}
  satisfy the simplifying internal consistency conditions
  \eqref{eq:consistency_conditions}.
\item[(ii)] If both the fast and slow methods have order at least two,
  then the overall multirate method is second order.
\item[(iii)] If both the fast and slow methods have order at least
  three, and if the outer method satisfies the additional condition
  \begin{align}
    \label{eq:rfsmr_3rd}
    \sum_{i=2}^{s^O} \(c_i^O-c_{i-1}^O\)
    \(\vect{e}_i+\vect{e}_{i-1}\)^{\intercal} A^O c^O
    + \(1-c_{s^O}^O\)
    \(\frac12+\vect{e}_{s^O}^{\intercal} A^O c^O\) = \frac13,
  \end{align}
  then the overall MIS is third order -- we note that this guarantees
  satisfaction of the third-order ``fast'' coupling condition,
  $
    \lamT{b}{f} \lam{A}{f,s} \lam{c}{s} = \frac16.
  $
\end{itemize}

\section{Relaxed Multirate Infinitesimal Step Methods}
\label{ch:RMIS}

In this work, we extend the MIS schemes to allow for methods with a
more `relaxed' definition of $\lam{b}{f}$.  All other
components remain identical to those in MIS methods, namely the
construction of $\lam{A}{f,f}$, $\lam{A}{f,s}$,
$\lam{A}{s,f}$, $\lam{A}{s,s}$ and $\lam{b}{s}$.  As with MIS methods,
we consider a pair of base methods, $T_I$ and $T_O$, corresponding
to the inner and outer base Runge-Kutta tables, having $s^I$
and $s^O$ stages, respectively, where $T_O$ is explicit and $T_I$ may
be either explicit or implicit; however, we require an additional
assumption that $T_I$ has an explicit first stage.
Due to this structural similarity between RMIS and MIS methods, we may
immediately leverage all aspects of the theory presented in
\cite{Gunther2013,Sandu2013multirate} that \emph{do not} utilize $\lam{b}{f}$:
\begin{itemize}
\item[(i)] The MIS and RMIS coefficients $\lam{A}{f,f}$, $\lam{A}{f,s}$,
  $\lam{A}{s,f}$ and $\lam{A}{s,s}$ satisfy the internal
  consistency conditions \eqref{eq:consistency_conditions}.
\item[(ii)] If both $T_I$ and $T_O$ have order at least two,
  then the RMIS method satisfies the ``slow'' second-order conditions
  (i.e.~\eqref{eq:mvOrderGARK1}-\eqref{eq:mvOrderGARK2} with $\sigma=s$).
\item[(iii)] If both $T_I$ and $T_O$ have order at least three, then
  the RMIS method satisfies the ``slow'' third-order conditions
  (i.e.~\eqref{eq:mvOrderGARK1}-\eqref{eq:mvOrderGARK3b} with
  $\sigma=s$ and $\nu=\{f,s\}$).
\end{itemize}
The overall order of accuracy for an RMIS multirate method
therefore depends both on the base methods and on the choice of
coefficients $\lam{b}{f}$.  To this end, we first extend the analysis
from G{\"u}nther and Sandu \cite{Gunther2013,Sandu2013multirate}
to consider the ``slow'' fourth-order conditions for MIS methods.

\begin{theorem}
\label{thm:RMIS_slow_order4}
Assume that the inner base method $T_I$ is at least third order, the
outer base method $T_O$ is at least fourth order, and that both
satisfy the row-sum consistency conditions
\eqref{eq:consistency_conditions}.  If $T_O$ is explicit and satisfies
the additional condition
\begin{equation}
\label{eq:rmis_4th}
  v^{O\intercal} A^O c^O = \frac{1}{12},
\end{equation}
where
\[
  v^O_i = \begin{cases}
    0,& i=1,\\
    b_i^O\(c_i^O - c_{i-1}^O\) + \(c_{i+1}^O - c_{i-1}^O\)\sum_{j=i+1}^{s^O}b_j^O,&1<i<s^O,\\
    b_{s^O}^O\(c_{s^O}^O-c_{s^O-1}^O\),&i=s^O,
  \end{cases}
\]
then the MIS and RMIS coefficients $\lam{A}{f,f}$, $\lam{A}{f,s}$,
$\lam{A}{s,f}$, $\lam{A}{s,s}$ and $\lam{b}{s}$ satisfy all of
the ``slow'' fourth-order conditions
(i.e.,~\eqref{eq:mvOrderGARK4a}-\eqref{eq:mvOrderGARK4d} with
$\sigma=s$ and $\nu,\mu=\{f,s\}$).
\end{theorem}

A detailed proof of this theorem is provided in \ref{sec:Proofs}.  However, we
note that this result, in combination with \cite[Theorem 3.1]{Sandu2013multirate},
guarantees that when using any third-order inner method
$T_I$, and any explicit fourth order outer method $T_O$ that satisfies
\eqref{eq:rmis_4th}, the corresponding MIS and RMIS methods automatically
satisfy all of the ``slow'' fourth-order order conditions.

We now turn our attention to the fast solution coefficients,
$\lam{b}{f} \in \R^{s^f}$.  Assuming that we select $T_O$ and $T_I$
according to the above criteria, then we must only select $\lam{b}{f}$
to satisfy the 14 ``fast'' order condition equations,
\begin{align}
  \lamT{b}{f}\ones^{\{s^f\}} &= 1, \label{eq:RMIS_fast_order1} \\
  \lamT{b}{f}\lam{c}{f} &= \tfrac12, \label{eq:RMIS_fast_order2} \\
  \lamT{b}{f}\(\lam{c}{f}\times\lam{c}{f}\) &= \tfrac13, \label{eq:RMIS_fast_order3a} \\
  \lamT{b}{f}\lam{A}{f,\nu}\lam{c}{\nu} &= \tfrac16, \label{eq:RMIS_fast_order3b} \\
  \lamT{b}{f}\(\lam{c}{f}\times\lam{c}{f}\times\lam{c}{f}\) &= \tfrac14, \label{eq:RMIS_fast_order4a} \\
  \(\lam{b}{f}\times\lam{c}{f}\)^{\intercal}\lam{A}{f,\nu}\lam{c}{\nu} &= \tfrac18, \label{eq:RMIS_fast_order4b} \\
  \lamT{b}{f}\lam{A}{f,\nu}\(\lam{c}{\nu}\times\lam{c}{\nu}\) &= \tfrac{1}{12}, \label{eq:RMIS_fast_order4c} \\
  \lamT{b}{f}\lam{A}{f,\mu}\lam{A}{\mu,\nu}\lam{c}{\nu} &= \tfrac{1}{24}, \label{eq:RMIS_fast_order4d}
\end{align}
where $\nu, \mu = \{f,s\}$, to guarantee that the overall RMIS method
is fourth-order accurate.  We note that each of these equations
depends \emph{linearly} on $\lam{b}{f}$, indicating that if
$s^f=s^Os^I$ is ``large enough'', we may be able to solve this linear
system for valid sets of coefficients $\lam{b}{f}$.  We further note
that for the RMIS method to be truly ``multirate'',
$s^f \propto m\,s^s$, indicating that as $m$ increases (i.e., the
problem has a larger multirate scale separation factor),
this linear system of equations becomes
severely under-determined, and hence the fourth-order conditions
become even simpler to satisfy.

While there exist multiple choices for $\lam{b}{f}$ that satisfy
the above criteria, we select
\begin{equation}
  \label{eq:RMIS_bf}
  \lamT{b}{f} = \begin{bmatrix}
    b_1^O \vect{e}_1^{\intercal} & b_2^O \vect{e}_1^{\intercal} &
    \cdots & b_{s^O}^O \vect{e}_1^{\intercal} \end{bmatrix} \in
  \R^{s^Os^I} = \R^{s^f}.
\end{equation}
There are three significant benefits of this choice: due to its
structured form this choice of $\lam{b}{f}$ may be used for arbitrary
multirate factors $m$, this choice of $\lam{b}{f}$ ensures that the resulting
multirate method conserves linear invariants, and this particular structure
results in a dramatic simplification of the order conditions
\eqref{eq:RMIS_fast_order1}-\eqref{eq:RMIS_fast_order4d}, as seen in
the following Lemma.

\begin{lemma}
\label{lem:RMIS_bf_simplifications}
Suppose that the coefficients $\lam{b}{f}$ are chosen as in equation
\eqref{eq:RMIS_bf}, and that the inner Butcher table $T_I$ has
explicit first stage (i.e.~the first entry of $c^I$ and the first row
of $A^I$ are identically zero).  Then the following identities hold:
\begin{align*}
  \lamT{b}{f}\(\lam{c}{f}\)^q &= \lamT{b}{s} \(\lam{c}{s}\)^q,
                                \;\forall q\ge 0,\\
  \lamT{b}{f} \lam{A}{f,f} &= \lamT{b}{s} \lam{A}{s,f},\\
  \lamT{b}{f} \lam{A}{f,s} &= \lamT{b}{s} \lam{A}{s,s},\\
  \(\lam{b}{f}\times\lam{c}{f}\)^{\intercal} \lam{A}{f,f} &= \(\lam{b}{s}\times\lam{c}{s}\)^{\intercal} \lam{A}{s,f},\\
  \(\lam{b}{f}\times\lam{c}{f}\)^{\intercal} \lam{A}{f,s} &= \(\lam{b}{s}\times\lam{c}{s}\)^{\intercal} \lam{A}{s,s},
\end{align*}
where $\lam{A}{f,f}$, $\lam{A}{s,f}$ and $\lam{A}{f,s}$ are defined as
in equations \eqref{eq:MIS_Aff}-\eqref{eq:MIS_Afs}, $\lam{c}{f}$ is
defined as in \eqref{eq:MIS_cf}, $\lam{c}{s} = c^O$, and $\lam{A}{s,s} = A^O$.
\end{lemma}

Again, we provide proof of this result in \ref{sec:Proofs}.

\begin{theorem}
\label{thm:RMIS_order4}
Let the base methods $T_I$ and $T_O$ satisfy all requirements for
Theorem \ref{thm:RMIS_slow_order4}, and assume that the coefficients
$\lam{b}{f}$ are chosen according to the formula \eqref{eq:RMIS_bf}.
If the first stage of $T_I$ is explicit, then the RMIS method defined via
\eqref{eq:MIS_Aff}-\eqref{eq:MIS_Afs} and \eqref{eq:RMIS_bf} is
fourth-order accurate. \\

The proof of this result follows directly from application of Theorem
\ref{thm:RMIS_slow_order4} and Lemma
\ref{lem:RMIS_bf_simplifications}.  From this lemma, all of the
``fast'' fourth-order conditions are equivalent to the corresponding
``slow'' order conditions:
\begin{align*}
  \lamT{b}{f}\lam{c}{f}\ &= \lamT{b}{s} \lam{c}{s} = \frac12,\\
  \lamT{b}{f}\(\lam{c}{f}\)^2 &= \lamT{b}{s} \(\lam{c}{s}\)^2 = \frac13,\\
  \lamT{b}{f}\(\lam{c}{f}\)^3 &= \lamT{b}{s} \(\lam{c}{s}\)^3 = \frac14,\\
  \lamT{b}{f}\lam{A}{f,\nu}\lam{c}{\nu} &=
    \lamT{b}{s}\lam{A}{s,\nu}\lam{c}{\nu} = \tfrac16 \\
  \lamT{b}{f}\lam{A}{f,\nu}\(\lam{c}{\nu}\times\lam{c}{\nu}\) &=
    \lamT{b}{s}\lam{A}{s,\nu}\(\lam{c}{\nu}\times\lam{c}{\nu}\) = \tfrac{1}{12}, \\
  \lamT{b}{f}\lam{A}{f,\mu}\lam{A}{\mu,\nu}\lam{c}{\nu} &=
    \lamT{b}{s}\lam{A}{s,\mu}\lam{A}{\mu,\nu}\lam{c}{\nu} = \tfrac{1}{24},\\
  \(\lam{b}{f}\times\lam{c}{f}\)^{\intercal}\lam{A}{f,\nu}\lam{c}{\nu} &=
    \(\lam{b}{s}\times\lam{c}{s}\)^{\intercal}\lam{A}{s,\nu}\lam{c}{\nu} = \tfrac18,
\end{align*}
where $\nu, \mu = \{f,s\}$.  Furthermore, all of these ``slow'' order
conditions are automatically met due to Theorem
\ref{thm:RMIS_slow_order4}. \qed
\end{theorem}

We point out that the above result no longer requires the condition
\eqref{eq:rfsmr_3rd} on $T_O$, that is typically required for MIS
methods to achieve third-order accuracy, due to the alternate
structure of $\lam{b}{f}$ within the RMIS approach.
We further note that based on the proof above, if all of the
assumptions from Theorem \ref{thm:RMIS_order4} are satisfied
\emph{except} for the condition \eqref{eq:rmis_4th} on $T_O$, then the
RMIS method will be third-order accurate.

\subsection{RMIS with MIS embedding}

Since RMIS methods and MIS methods share the same algorithmic
structure, have identical coefficients $\lam{A}{f,f}$, $\lam{A}{f,s}$,
$\lam{A}{s,f}$, $\lam{A}{s,s}=A^O$ and $\lam{b}{s}=b^O$, and only
differ in their selection of $\lam{b}{f}$, one may naturally question
whether an MIS method could be used as an embedding within a RMIS
method to enable temporal error estimation.  For this combination
of RMIS and MIS methods to work, however, the base methods $T_O$ and
$T_I$ must be compatible with both approaches.  Specifically, if
these base methods satisfy the requirements:
\begin{itemize}
\item[(a)] $T_O$ is explicit and at least order four,
\item[(b)] $T_I$ has explicit first stage and is at least order three,
\item[(c)] $T_O$ satisfies the MIS condition \eqref{eq:rfsmr_3rd},
  i.e.
  \[
    \sum_{i=2}^{s^O} \(c_i^O-c_{i-1}^O\) \(\vect{e}_i+\vect{e}_{i-1}\)^{\intercal} A^O c^O
    + \(1-c_{s^O}^O\) \(\frac12+\vect{e}_{s^O}^{\intercal} A^O c^O\) = \frac13,
  \]
\item[(d)] $T_O$ satisfies the RMIS condition \eqref{eq:rmis_4th},
  i.e. $v^{O\intercal} A^O c^O = \frac{1}{12}$, where
  \[
    v^O_i = \begin{cases}
      0,& i=1,\\
      b_i^O\(c_i^O - c_{i-1}^O\) + \(c_{i+1}^O - c_{i-1}^O\)\sum_{j=i+1}^{s^O}b_j^O,&1<i<s^O,\\
      b_{s^O}^O\(c_{s^O}^O-c_{s^O-1}^O\),&i=s^O,
    \end{cases}
  \]
\end{itemize}
then the combination of shared MIS+RMIS coefficients $\lam{A}{f,f}$,
$\lam{A}{f,s}$, $\lam{A}{s,f}$, $\lam{A}{s,s}$ and $\lam{b}{s}$, along
with the RMIS coefficients $\lam{b}{f}$ given by equation
\eqref{eq:RMIS_bf} and MIS coefficients
\begin{equation}
  \label{eq:RMIS_btildef}
  \tilde{\vect{b}}^{\{f\}\intercal} = \begin{bmatrix}
    c_2^O b^{I\intercal} & \(c_3^O-c_2^O\) b^{I\intercal} & \cdots &
    \(1-c_{s^O}^O\)b^{I\intercal} \end{bmatrix} \in \R^{s^f},
\end{equation}
will result in a fourth-order method with third-order embedding.

This embedding allows temporal error estimation for local truncation
error without requiring extra function evaluations.  Since both the
MIS and RMIS method solutions may be constructed using the same stage
values, their difference provides an $\mathcal{O}(h^3)$-accurate
estimate of the overall multirate solution error. Based on this error
estimate, standard temporal adaptivity controllers may be used to
accept/reject the current `slow' time step, as well as estimate an
`optimal' step size $h$ for the ensuing step.  We note, however, that
since this is only an estimate of the overall multirate solution
error, it cannot be used to simultaneously adapt the time scale
separation factor $m$, thus alternate mechanisms for dual adaptivity
of both $h$ and $m$ are an ongoing area of research.

\begin{figure}[tbph]
  \centerline{
    \includegraphics[width=.6\linewidth]{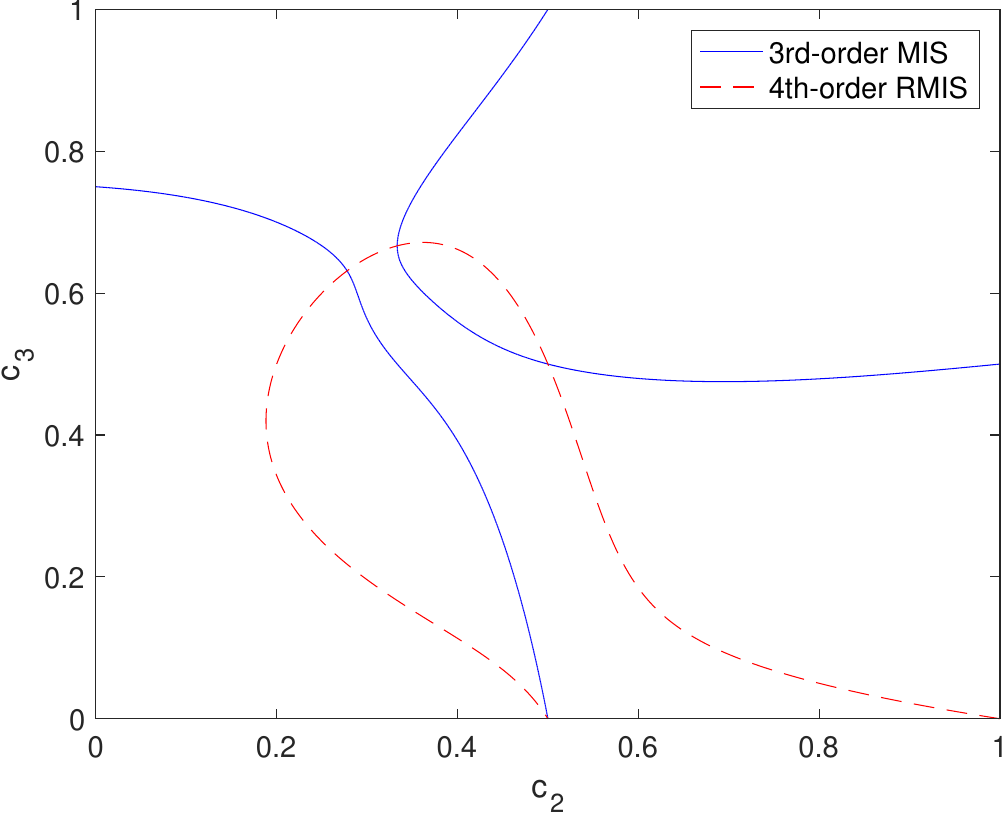}
    }
    \caption{Choices of $c_2$ and $c_3$ that satisfy the conditions
      \eqref{eq:RK4stageSolRFSMR} and \eqref{eq:RK4stageSolRMIS}.
      Values of $c_2$ and $c_3$ \emph{not} on these curves result in
      MIS and RMIS methods of 2nd and 3rd order, respectively.
    }
    \label{fig:RK4stageSol}
\end{figure}

While the above criteria may be met by a variety of base methods, we
identify some specific candidates that we use for our numerical
results in Section \ref{sec:numerical_tests}.  To this end, we relied
on Butcher's derivation of families of explicit 4th order methods
\cite{Butcher2008} to determine $T_O$.  His general solution for a
4-stage fourth order explicit Runge-Kutta method depends on only two
free variables, $(c_2,c_3)$, and is reproduced in
\ref{sec:Coefficients}.  For $T_O$ with this structure, the 3rd-order
MIS condition \eqref{eq:rfsmr_3rd} simplifies to
\eqref{eq:RK4stageSolRFSMR}, while the 4th-order RMIS condition
\eqref{eq:rmis_4th} simplifies to \eqref{eq:RK4stageSolRMIS}.
We plot the solutions $(c_2,c_3)\in[0,1]^2$ to these equations in
Figure \ref{fig:RK4stageSol}, and point out that the intersection of
these curves denote choices of $T_O$ that satisfy \emph{all} of our
desired criteria (a)-(d).
Of these, only two satisfy the MIS criteria that $c_2\le c_3$:
\begin{align}
  \label{eq:cvals1}
  (c_2,c_3) &= \(\frac13, \frac23\), \\
  \label{eq:cvals2}
  (c_2,c_3) &= \(\frac{2502984374488603}{9007199254740992}, \frac{2843567935040037}{4503599627370496}\)\\
  \notag
  &\approx \(0.27788708828342423285, 0.63139891871345210639\).
\end{align}
The first of these corresponds to the ``3/8-Rule'' from Kutta's 1901 paper
\cite{Kutta1901},
\begin{equation}
  \label{eq:RK4stage3/8}
  \begin{array}{c|cccc}
    0&    0&   0&   0&   0\\
    \frac13 &  \frac13 &   0&   0&   0\\
    \frac23 & -\frac13 &   1&   0&   0\\
    1&    1&  -1&   1&   0\\
    \hline
    &  \frac18& \frac38& \frac38& \frac18
  \end{array}
\end{equation}
Both due to its historical elegance, as well as its equally-spaced
abcissae $c_i$, we use this base method for $T_O$ in our numerical
results.

\subsection{RMIS example methods}
For numerical testing, we examine methods wherein the inner table
$T_I$ corresponds to a subcycled version (using approximately $m/s^O$
substeps) of the outer table $T_O$.  This ``telescopic'' approach
facilitates further recursion to support problems with three
or more rates; however, that is not studied in this work.  We note
that due to the patterned structure of both the MIS and RMIS methods,
these subcycled algorithms may be implemented more efficiently than a
generic GARK scheme having $ms$ fast and $s$ slow stages.

For these methods, we utilize $T_O$ as either the four-stage 3/8-Rule
\eqref{eq:RK4stage3/8} above, or the three-stage ``KW3'' table,
\begin{equation}
   \label{eq:RK3stageKW3}
   \begin{array}{c|ccc}
     0 & 0 & 0 & 0 \\
     \frac13 & \frac13 & 0 & 0 \\
     \frac34 & -\frac{3}{16} & \frac{15}{16} & 0 \\
     \hline
             & \frac16 & \frac{3}{10} & \frac{8}{15} \\
   \end{array}
\end{equation}
that satisfies the 3rd-order MIS condition \eqref{eq:rfsmr_3rd}.  We
note that this latter table is frequently used for MIS
methods in the literature
\cite{Schlegel2012,Knoth2014,Schlegel2009,Wensch2009}; moreover, in
previous tests of a wide range of multirate methods, we found this
combination to be the most efficient (smallest error with least
computational cost).

To be more precise regarding how we form $T_I$ from these two choices
of $T_O$ to enable subcycling, we elaborate on the case of a time-scale
separation factor $m=100$.  Schlegel suggested using
$n_{i}=\left\lceil m\left(c^O_{i}-c^O_{i+1}\right)\right\rceil $ to measure
the appropriate number of subcycles for each subcycling period
\cite{Schlegel2012b}.
We used this measure as a guide when selecting $m$ in our methods,
although we determined that our method comparisons would be fairest if
the ratio of fast function evaluations to slow function evaluations
was held approximately constant.

The 3/8-Rule \eqref{eq:RK4stage3/8} has four evenly-spaced abcissae
with final $c^O_{s^O}=1$.  Hence its multirate implementation requires only
three subcycling periods, and so we form $T_I$ using 34 subcycles of
the $T_O$ table.  The KW3 method has only three abcissae, but
$c^O_{s^O}=\frac34 < 1$, so although it has one fewer stage than the
3/8-Rule, it requires the same number of subcycling periods.
Moreover, since these abcissae are unevenly spaced, we could
choose between using either different $A^I$ per outer stage (to attempt
nearly-identical substep sizes) or using a single $A^I$ for each outer
stage, mapped to the largest subcycling interval.  Choosing the latter
approach for simplicity, we form $T_I$ using 35 subcycles of
$T_O$. This is a compromise between using the largest subcycling
interval for substep sizes and making the work of fast functions vs
slow functions comparable.

This gives us the multirate methods:
\begin{itemize}
\item ``RMIS-3/8'': this is our proposed RMIS method using
  $T_O$ given by \eqref{eq:RK4stage3/8} that satisfies both the
  4th-order condition \eqref{eq:rmis_4th} and the 3rd-order condition
  \eqref{eq:rfsmr_3rd} -- the method should be $\mathcal O(h^4)$ accurate;
\item ``RMIS-KW3'': this is our proposed RMIS method using
  $T_O$ given by \eqref{eq:RK3stageKW3} that satisfies
  \eqref{eq:rfsmr_3rd} but does not satisfy \eqref{eq:rmis_4th} -- the
  method should be $\mathcal O(h^3)$ accurate.
\end{itemize}

\subsection{Optimized method}
\label{sec:opt}

To further explore the preceding theoretical results, we also consider
an ``unstructured'' form of $\lam{b}{f}\in\R^{s^f}$, where these
coefficients are chosen to directly satisfy the equations
\eqref{eq:RMIS_fast_order1}-\eqref{eq:RMIS_fast_order4d}.  Here, since
we do not enforce any particular structure on these coefficients, we
must consider a fixed time scale separation, $m$ (or equivalently
$s^f$), here chosen to be $m=100$.  As with the RMIS methods above, we
choose $T_I$ as a subcycled version of $T_O$, which we again choose to
be the 3/8-Rule, resulting in $s^f = 408$.  Since this allows
significantly more coefficients than the number of fourth order
condition equations, we use the remaining degrees of freedom in this
under-determined system to minimize the two-norm of the residual of all
the fifth-order conditions.  To solve this minimization problem, we used
MATLAB's \texttt{fminsearch} algorithm, which is based on a simplex search
method \cite{Lagarias1998}.  Under this approach, we arrive at the
multirate method ``Opt-3/8,'' which should be $\mathcal O(h^4)$
accurate.

\subsection{Method Comparison}

In the sections that follow, we examine the linear stability, and
assess the attainable order of accuracy, for both the proposed RMIS
and MIS algorithms.
In the ensuing sections we therefore compare the 3 newly developed
multirate methods listed above with 2 existing multirate methods:
\begin{itemize}
\item ``MIS-3/8'': this is the MIS method using
  $T_O$ given by \eqref{eq:RK4stage3/8} -- the resulting method should
  be $\mathcal O(h^3)$ accurate;
\item ``MIS-KW3'': this is the MIS method using
  $T_O$ given by \eqref{eq:RK3stageKW3} -- the resulting method should
  be $\mathcal O(h^3)$ accurate.
\end{itemize}

We note that although we have suggested an embedding of MIS within
RMIS (i.e., MIS-3/8 within RMIS-3/8), optimal approaches for time-step
adaptivity (that vary both $h$ and $m$) are not within the scope of
this manuscript and are left for future work.  We therefore explore
each of the above methods using fixed time-step sizes $h$.

\section{RMIS Implementation and Memory Considerations}
\label{ssec:Memory}

A crucial aspect for the success of multirate methods is their ability
to be implemented in both a computationally and memory-efficient
manner.  First, we recall the formulas \eqref{eq:GARK_step},
corresponding to the generic implementation of a 2-component GARK
method:
\begin{align*}
  \slamc{k}{f}{j} &= y_{n} +
                     h \sum_{l=1}^{s^{\{f\}}} \slamc{a}{f,f}{j,l} f^{\{f\}}\(\slamc{k}{f}{l}\) +
                     h \sum_{l=1}^{s^{\{s\}}} \slamc{a}{f,s}{j,l} f^{\{s\}}\(\slamc{k}{s}{l}\), \\
  \slamc{k}{s}{i} &= \vect{y}_{n} +
                     h \sum_{l=1}^{s^{\{f\}}} \slamc{a}{s,f}{i,l} f^{\{f\}}\(\slamc{k}{f}{l}\) +
                     h \sum_{l=1}^{s^{\{s\}}} \slamc{a}{s,s}{i,l} f^{\{s\}}\(\slamc{k}{s}{l}\),\\
  y_{n+1} &= y_{n} + h \sum_{l=1}^{s^{\{f\}}} \slamc{b}{f}{l} f^{\{f\}}\(\slamc{k}{f}{l}\) +
                    h \sum_{l=1}^{s^{\{s\}}} \slamc{b}{s}{l} f^{\{s\}}\(\slamc{k}{s}{l}\),
\end{align*}
where $i=1,\dots,s^{\{s\}}$ and $j=1,\dots,s^{\{f\}}$.
We introduce additional notation to better focus on the case where
subcycling is employed in the inner method $T_I$.  Subcycling suggests
a natural set of internal divided units, that can be characterized by
the block-rows of $\lam{A}{f,f}$ and $\lam{A}{f,s}$ which correspond
to each fast substep.  To this end, we define the fast stage abcissae
in vector form as
\[
  \lam{c}{f}=\begin{bmatrix}\lam{c}{f,1}\\
    \vdots\\
    \lam{c}{f,s^{O}}
  \end{bmatrix},
  \quad \text{where} \quad
  \lam{c}{f,i}=\begin{bmatrix}c_{i-1}^{O}+c_{1}^{I}\left(c_{i}^{O}-c_{i-1}^{O}\right)\\
    \vdots\\
    c_{i-1}^{O}+c_{s^{I}}^{I}\left(c_{i}^{O}-c_{i-1}^{O}\right)
  \end{bmatrix}.
\]
We make row block indices corresponding to these stage vectors, which
correspond to the smallest unit step; the column block indices
delineate the same organization for how $\lam{b}{f}$ is used.
This then gives us $\lam{A}{f,f,i,j}$ and $\lam{A}{f,s,i}$, where
$i$ denotes the block-row, and $j$ denotes the block-column. This
also allows $\slam{k}{f}$ to be represented in terms of block-row
units,
$\slam{k}{f}={\begin{bmatrix}{\slam{k}{f,1}}^{\intercal} & \cdots &
    {\slam{k}{f,s^{O}}}^{\intercal} \end{bmatrix}}^{\intercal}$,
where $\slamc{k}{f,i}{j}$ is the $j$th fast stage vector for the $i$th block.
Since we have assumed the first stage of $T_I$ is explicit, then the
coefficients $\slamc{a}{f,f,i,j}{1,q}=\slamc{a}{s,f,j}{i,q}$ and
$a_{1,l}^{\{f,s,i\}}=a_{i,l}^{\{s,s\}}$.
We leverage this structure by specifically defining
$\slamc{\hat{k}}{f,i}{1}$ (the fast stage solution from the $i$th
block) to be $\slamc{k}{f,i}{1}$ with the slow coupling portion
removed,
\[
  \slamc{\hat{k}}{f,i}{1} = \slamc{k}{f,i}{1} -
  h\sum_{l=1}^{s^{\{s\}}}a_{1,l}^{\{f,s,i\}} f^{\{s\}}\(\slamc{k}{s}{l}\).
\]
This portion of the stage solution is then used to accumulate an
initial condition for $\slamc{k}{f,i}{p}$ for $p=2,\dots,s^I$,
and to calculate the next initial condition
$\slamc{\hat{k}}{f,i+1}{1}$.  Finally, this accumulated initial
condition $\slamc{k}{f,s^O}{1}$ can be used with the remaining
$\slamc{k}{f,s^O}{p}$ stages to form the fast portion of an embedded
MIS solution method (if desired).  The formulas for these are as
follows:

\small
\begin{align*}
  \slamc{\hat{k}}{f,1}{1} &=  y_{n} \\
  \slamc{k}{f,i}{1}       & =  \slamc{\hat{k}}{f,i}{1} + h\sum_{l=1}^{s^{\{s\}}}a_{1,l}^{\{f,s,i\}} f^{\{s\}}\(\slamc{k}{s}{l}\)\\
  \slamc{\hat{k}}{f,i}{1} & =  \sum_{p=1}^{i-1}\slamc{k}{f,i}{1} +
                           h\sum_{l=1}^{s^I}b_{l}^I\(c_{i}^O-c_{i-1}^O\)
                           f^{\{f\}}\(\slamc{k}{f,i-1}{l}\)\\
  \slamc{k}{f,i}{j} & = \slamc{\hat{k}}{f,i}{1} + h\sum_{l=1}^{j-1}a_{jl}^{\{f,f,i,i\}} f^{\{f\}}\(\slamc{k}{f}{l}\) +
                     h\sum_{l=1}^{s^{\{s\}}}a_{jl}^{\{f,s,i\}} f^{\{s\}}\(\slamc{k}{s}{l}\)\\
  \slamc{k}{s}{i} & = \slamc{\hat{k}}{f,i}{1} + h\sum_{l=1}^{i-1}a_{il}^{\{s,s\}} f^{\{s\}}\(\slamc{k}{s}{l}\)\\
  \tilde{y}_{n+1}&=\slamc{k}{f,s^O}{1} + h\sum_{l=1}^{s^I} b_{l}^I\(1-c_{s^O}^O\) f^{\{f\}}\(\slamc{k}{f,s^O}{l}\)
                 + h\sum_{l=1}^{s^{\{s\}}} b_l^{\{s\}} f^{\{s\}}\(\slamc{k}{s}{l}\)\\
  y_{n+1} & = y_n + h\sum_{l=1}^{s^{\{f\}}} b_l^{\{f\}} f^{\{f\}}\(\slamc{k}{f}{l}\) +
                   h\sum_{l=1}^{s^{\{s\}}} b_l^{\{s\}} f^{\{s\}}\(\slamc{k}{s}{l}\)
\end{align*}
\normalsize
where the final step solution $y_{n+1}$ is built up by successive fast
substeps, and where the slow information is updated as needed.

To remove duplicate function calls, we use temporary vectors for the
contributions of our base methods to the solution.  In particular,
while computing the update formulas for
$\slamc{k}{f,i}{j}$ and $\slamc{k}{s}{i}$, we store $\slamc{k}{f,i}{1}$,
the $s^I$ fast function vectors $\slamc{\phi}{f,i}{l}=f^{\{f\}}\(\slamc{k}{f,i}{l}\)$,
and the $s^O$ slow function vectors $\slamc{\phi}{s}{l}=f^{\{s\}}\(\slamc{k}{s}{l}\)$.
These provide sufficient storage to compute the embedded MIS solution
$\tilde{y}_{n+1}$ without additional function calls.  Moreover,
to compute the RMIS solution $y_{n+1}$, we must additionally
store $y_{n+1}$ and the vector of coefficients $\lam{b}{f}$, if
these are unique or unstructured.  We may condense storage slightly so
that we need only retain $s^I+s^O+2$ vectors
the size of our solution $y$, by overwriting $\slam{\phi}{f,i}$ with
$\slam{\phi}{f,i+1}$ after the $i$th fast substep has completed, and
$\slam{\phi}{f,i}$ has been used to generate the initial condition
$\slamc{\hat{k}}{f,i+1}{1}$ for the fast substep $(i+1)$.  When using
subcycling for the fast steps, we overwrite $\slam{\phi}{f,i}$ after each
subcycled fast substep has completed, and $\slam{\phi}{f,i}$ has been used
to update the temporary initial condition $\slamc{\hat{k}}{f,i}{1}$ for
the next subcycled fast substep.

When considering only the structured form of MIS and RMIS methods,
further optimizations may be used since the overall time step may be
broken apart into a sequence of ``fast'' subproblems, and where each
subproblem only receives forcing terms from the slow stage right-hand
side functions $f^{\{s\}}$.  To this end, we have provided an open-source
MATLAB implementation of this structured form for both RMIS and MIS
methods \cite{RMIS_repo}.  We note that this repository does not
include an implementation of the ``Opt-3/8'' method described
previously in Section \ref{sec:opt}, as that requires an unstructured
$\lam{b}{f}$; however, in this repository we also include a file
\texttt{opt38\_coeffs.m} that contains these $\lam{b}{f}$ coefficients
for interested readers.

\section{RMIS and MIS linear stability analysis}
\label{sec:linear_stability}

In their paper introducing GARK methods \cite{Sandu2015}, Sandu and
G{\"u}nther analyze linear stability using a modification of the
standard Dahlquist test problem,
\[
  y'= \sum_{m=1}^N \lambda^{\{m\}} y, \qquad y(0) = 1,
\]
which in the current context of a two-rate problem becomes
\[
  y' = \lambda^{\{f\}}y + \lambda^{\{s\}}y, \qquad y(0) = 1,
\]
and where they assume that the real parts of both $\lambda^{\{f\}}$
and $\lambda^{\{s\}}$ must be negative.  While that approach yields
elegant definitions of the resulting GARK stability regions, we prefer
an approach for multirate linear stability analysis that was first
proposed by Kv{\ae}rn{\o} \cite{Kvaerno2000} in 2000, has
subsequently been used by a variety of authors
\cite{Constantinescu2013,Savcenco2007,Savcenco2008,Hundsdorfer2009,Kuhn2014},
and which we summarize here.  In this approach, one instead considers
the partitioned test problem
\begin{align}
  \label{eq:partitioned_linear}
  &\begin{bmatrix} y_f' \\ y_s' \end{bmatrix} =
  \begin{bmatrix} g_{11} & g_{12}\\ g_{21} & g_{22} \end{bmatrix}
  \begin{bmatrix} y_f \\ y_s \end{bmatrix}, \qquad
  \begin{bmatrix} y_f(0) \\ y_s(0) \end{bmatrix} =
  \begin{bmatrix} 1 \\ 1 \end{bmatrix} \\
  \notag
  \Leftrightarrow\quad&\\
  \notag
  &y' = Gy, \qquad y(0) = \ones^{\{2\}}.
\end{align}
Here, $y_f$ and $y_s$ correspond to the fast and slow variables,
respectively, with a resulting splitting of the right-hand side into
fast and slow components as shown in
\eqref{eq:partitioned_reformulation}.  The benefit of this approach is
that unlike the purely additive scalar version, this directly allows
analysis of stability as a function of the strength of
fast/slow coupling in the problem, and does not require simultaneous
diagonalizability of the Jacobians of $\flam{s}$ and $\flam{f}$ to
derive the linear stability problem.  Defining the matrix $Z=hG$,
stability is ascertained by first determining the amplification
matrix $S(Z)$ for one step of the numerical method, i.e.
\[
  y_{n+1} = S(Z)\, y_n, \quad \text{where} \quad
  S(Z) = \begin{bmatrix} s_{11}(Z) & s_{12}(Z)\\ s_{21}(Z) &
    s_{22}(Z) \end{bmatrix}.
\]
The multirate method is then linearly stable if the eigenvalues
of $S$ have magnitude less than $1$ for a given step size $Z=hG$.

While more complex than standard IVP linear stability regions, the
number of independent parameters may be reduced from five
($h, g_{1,1}, g_{1,2}, g_{2,1}, g_{2,2}$) to a more
manageable three \cite{Savcenco2007,Savcenco2008,Kuhn2014}.  First,
the eigenvalues of $G$ have negative real part if $g_{11},g_{22}<0$ and
$\beta=\frac{g_{12}g_{21}}{g_{11}g_{22}}<1$, where here $\beta$ is a
measure of the off-diagonal coupling strength of the problem.
Under this assumption, the eigenvalues of $S(Z)$ depend on only
three parameters:
\begin{align}
  \kappa &= \frac{g_{22}}{g_{11}} > 0,\\
  \xi &= \frac{hg_{11}}{1-hg_{11}} \in \left(-1,0\right),\\
  \eta &= \frac{\beta}{2-\beta} \in (-1,1).
\end{align}
Here, $\kappa$ encodes the time-scale separation of the problem, $\xi$
encodes the stiffness of the fast time scale (as $hg_{11}$ goes from
$0\to-\infty$, $\xi$ moves from $0\to -1$), and $\eta$ encodes the
strength of coupling (no coupling at $\eta=0$, increased coupling as
$|\eta|\to 1$, coupling-dominant as $\eta\to -1$).  The stability of a
multirate method may then be visualized using snapshots of the
$(\xi,\eta)$ stability regions for fixed values of $\kappa$.

In the stability region plots that follow, we match the problem and
method time-scale separation values, i.e.~$\kappa = m$.  We then
create linearly-spaced arrays of 100 $\xi$ values in $(-1,0)$ and 200 $\eta$
values in $(-1,1)$.  For each of the 20000 resulting $(\xi_i,\eta_j)$
combinations, we compute the eigenvalues of $S(Z)$ to determine
stability, and mark the stable regions in yellow and unstable
regions in blue.

\begin{figure}[tbph!]
  \centerline{
    \includegraphics [width=.49 \textwidth]{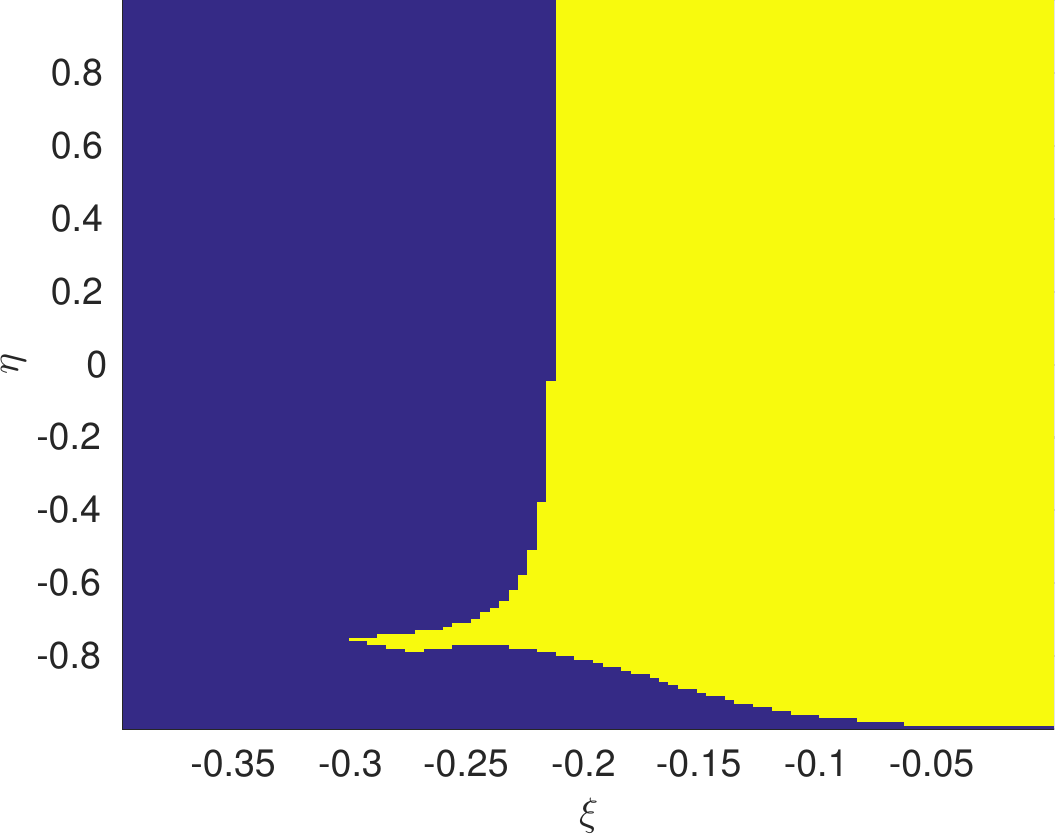}
    \hfill
    \includegraphics [width=.49 \textwidth]{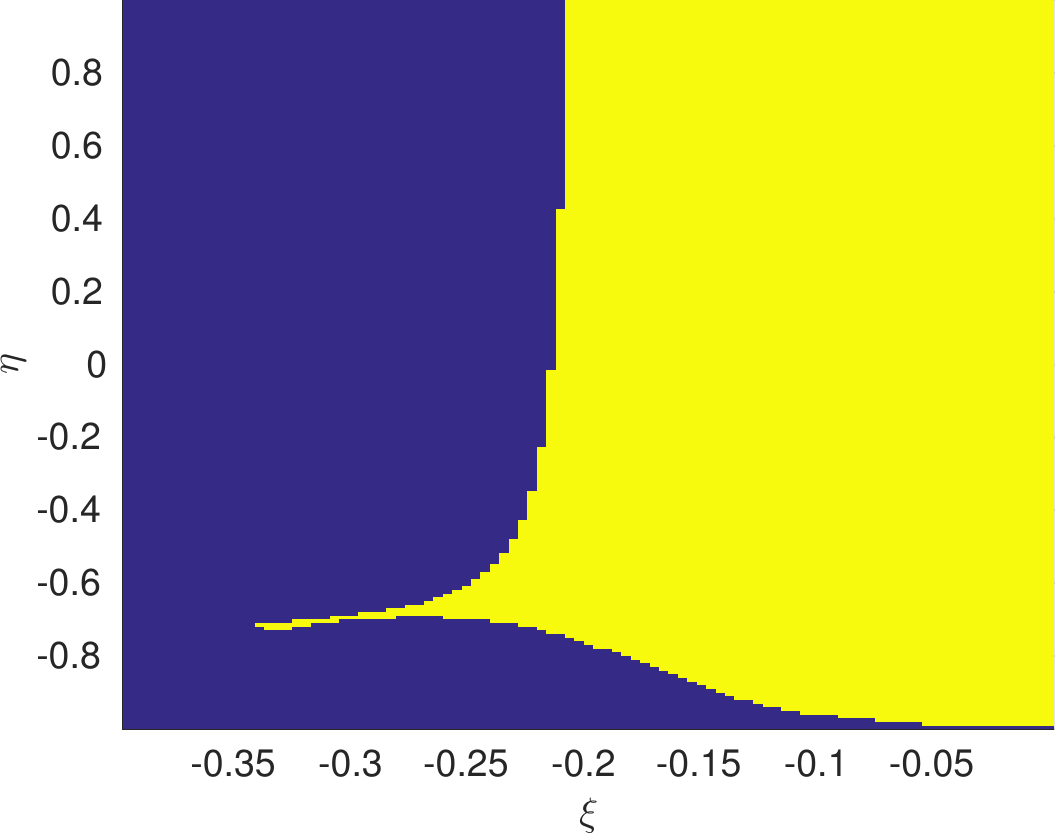}
  }
  \vspace{0.5em}
  \centerline{
    \includegraphics [width=.48 \textwidth]{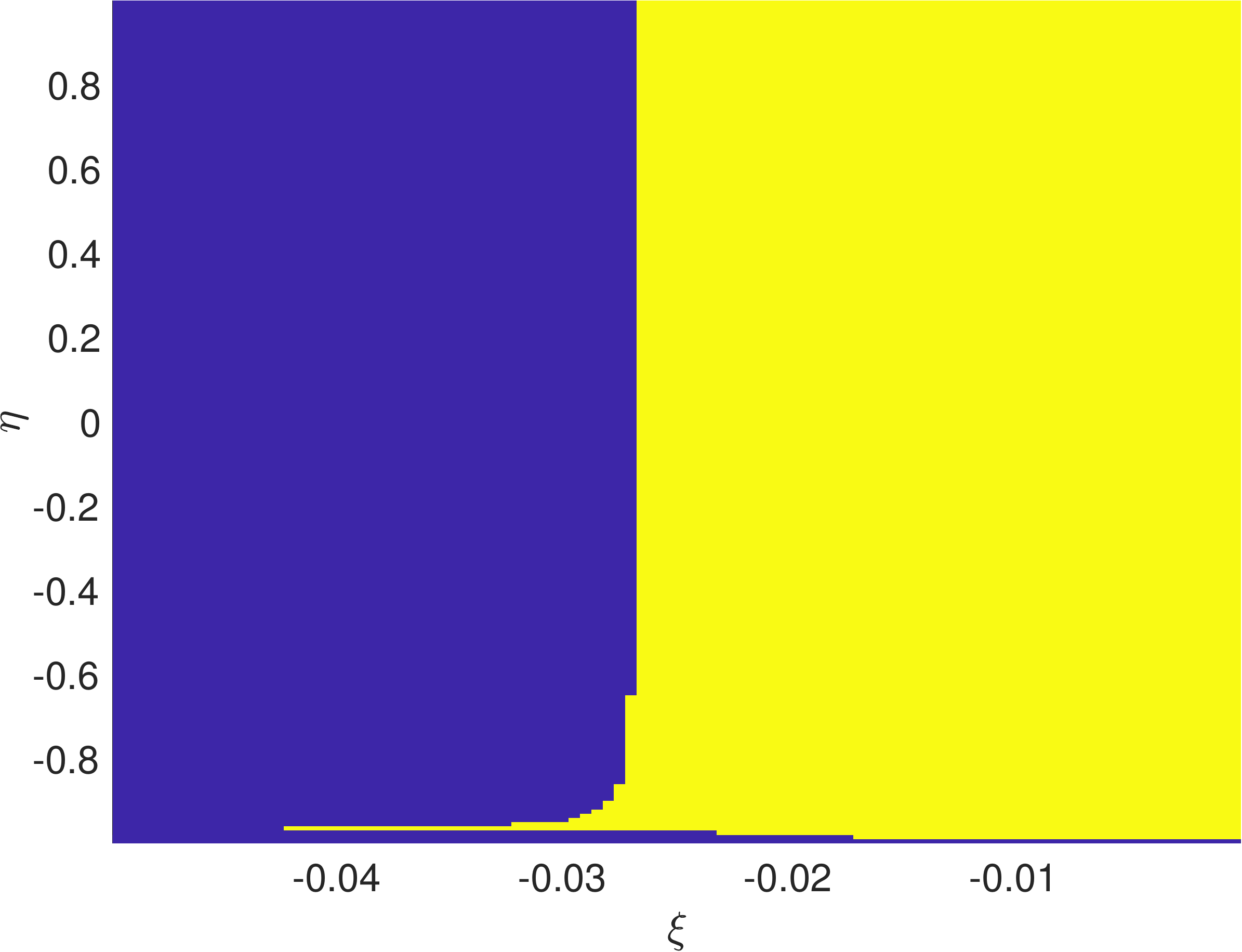}
    \hfill
    \includegraphics [width=.48 \textwidth]{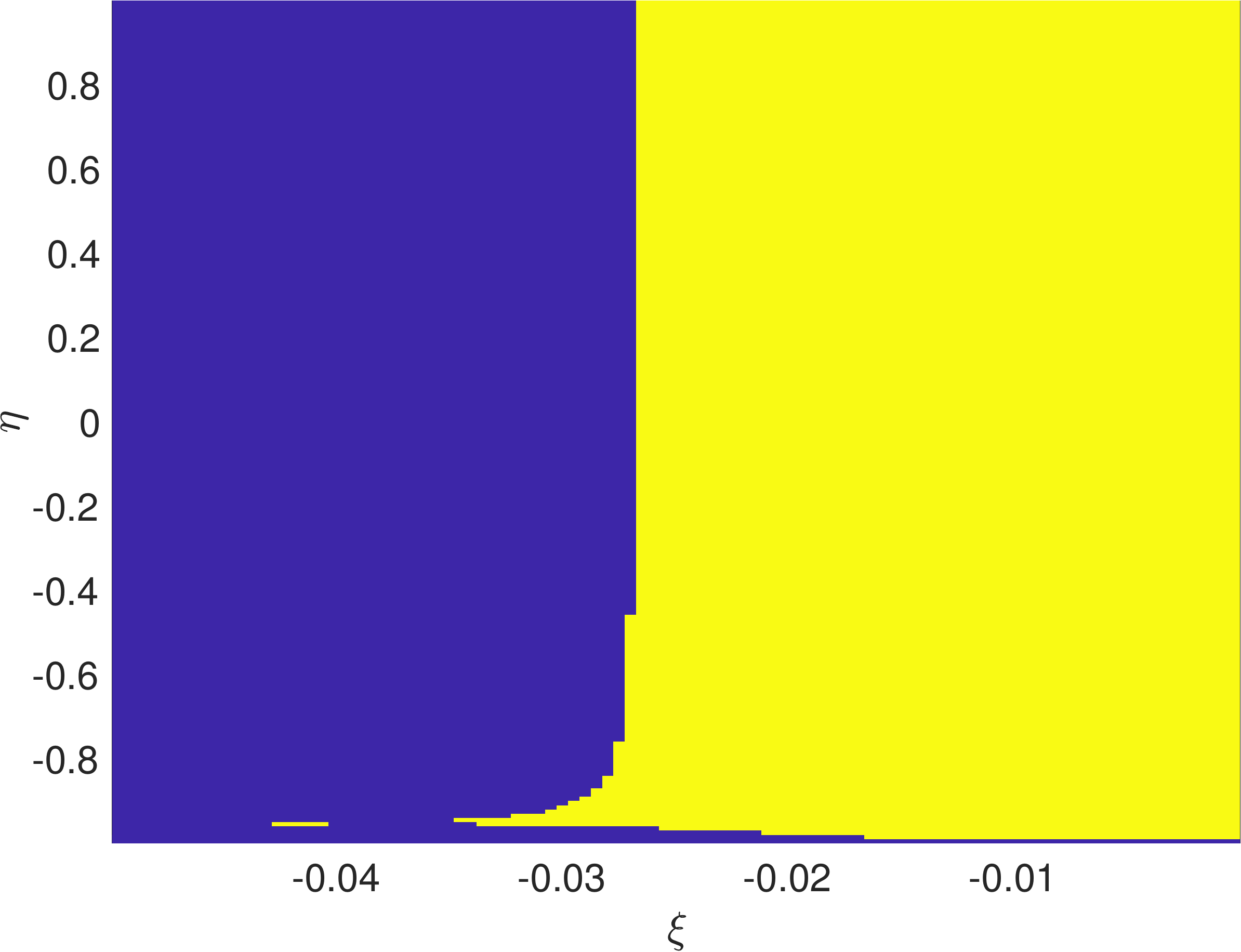}
  }
  \caption{Linear stability plots for the RMIS-3/8 (left) and
    MIS-3/8 (right) methods at time-scale separation values
    $\kappa=m=10$ (top) and $\kappa=m=100$ (bottom), showing the
    stable region in yellow.  Note that the $\kappa=100$ plots are
    zoomed in to a smaller $\xi$ region than the $\kappa=10$ plots.}
  \label{fig:case1stab}
\end{figure}

In Figure \ref{fig:case1stab} we plot the stability regions for the
RMIS-3/8 and MIS-3/8 methods for $\kappa=m=\{10,100\}$.  We see that
the stability regions for both methods are similar: for weakly coupled fast
and slow time scales, these indicate stability for $\xi$ values close to
zero and instability for $\xi$ close to -1, as would be expected for any
explicit method.  Additionally, both exhibit a region of increased
stability for problems with stronger fast/slow coupling, followed by
complete instability as $\eta\to -1$ (where the coupling terms become
infinitely large).  Comparing the plots of $\kappa=10$ versus 100, we
see that the region of increased stability shifts toward stronger
fast/slow coupling as the problem time-scale separation increases.
Lastly, as expected for explicit methods, the maximum stable step size
shrinks as the fast time scale decreases.

\begin{figure}[tbph!]
  \centerline{
    \includegraphics [width=.49 \textwidth]{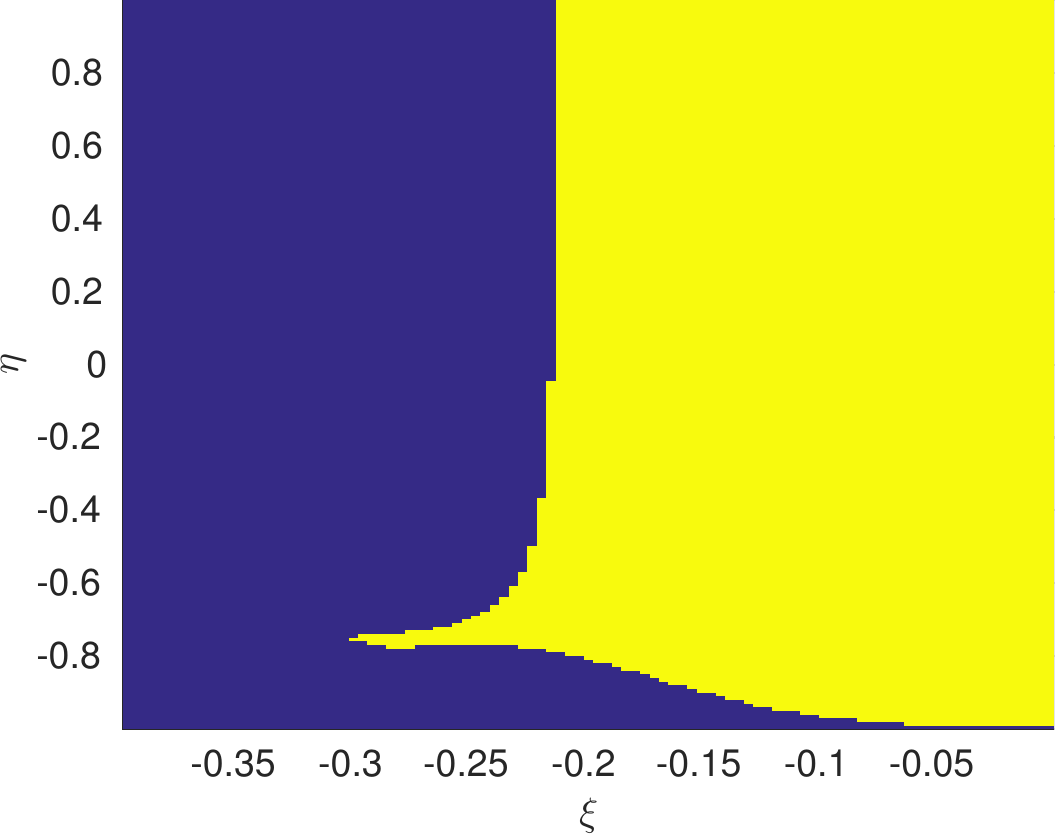}
    \hfill
    \includegraphics [width=.49 \textwidth]{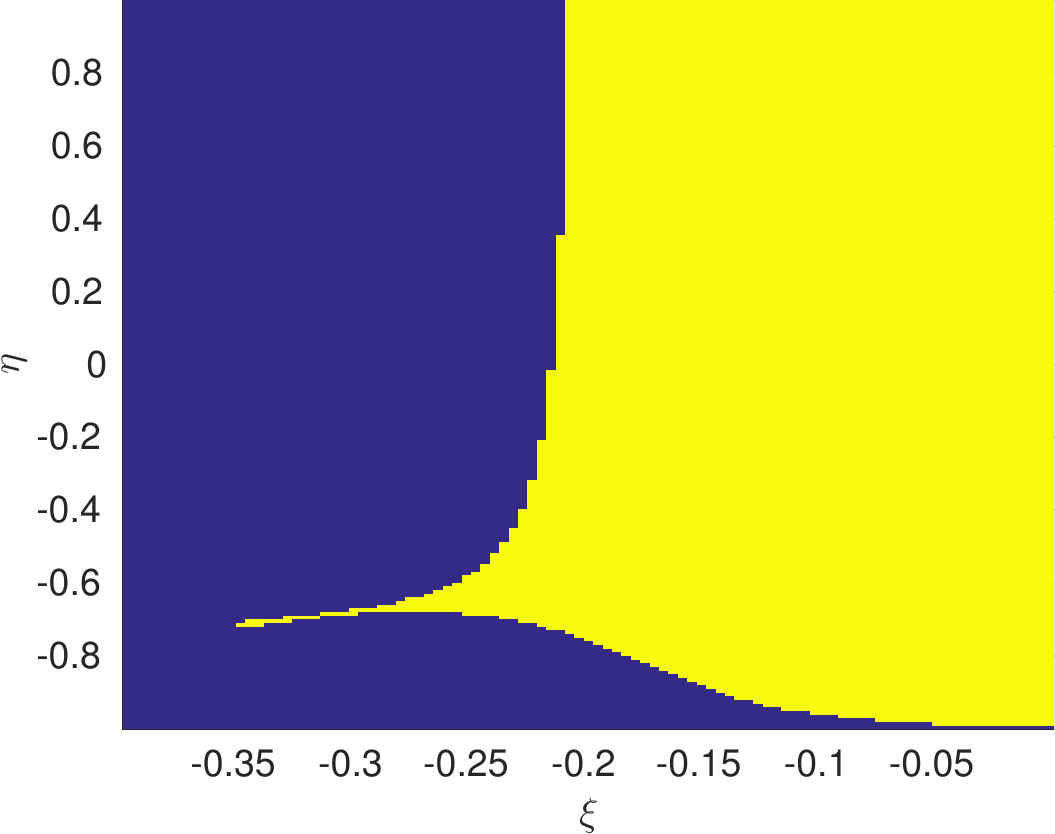}
  }
  \vspace{0.5em}
  \centerline{
    \includegraphics [width=.48 \textwidth]{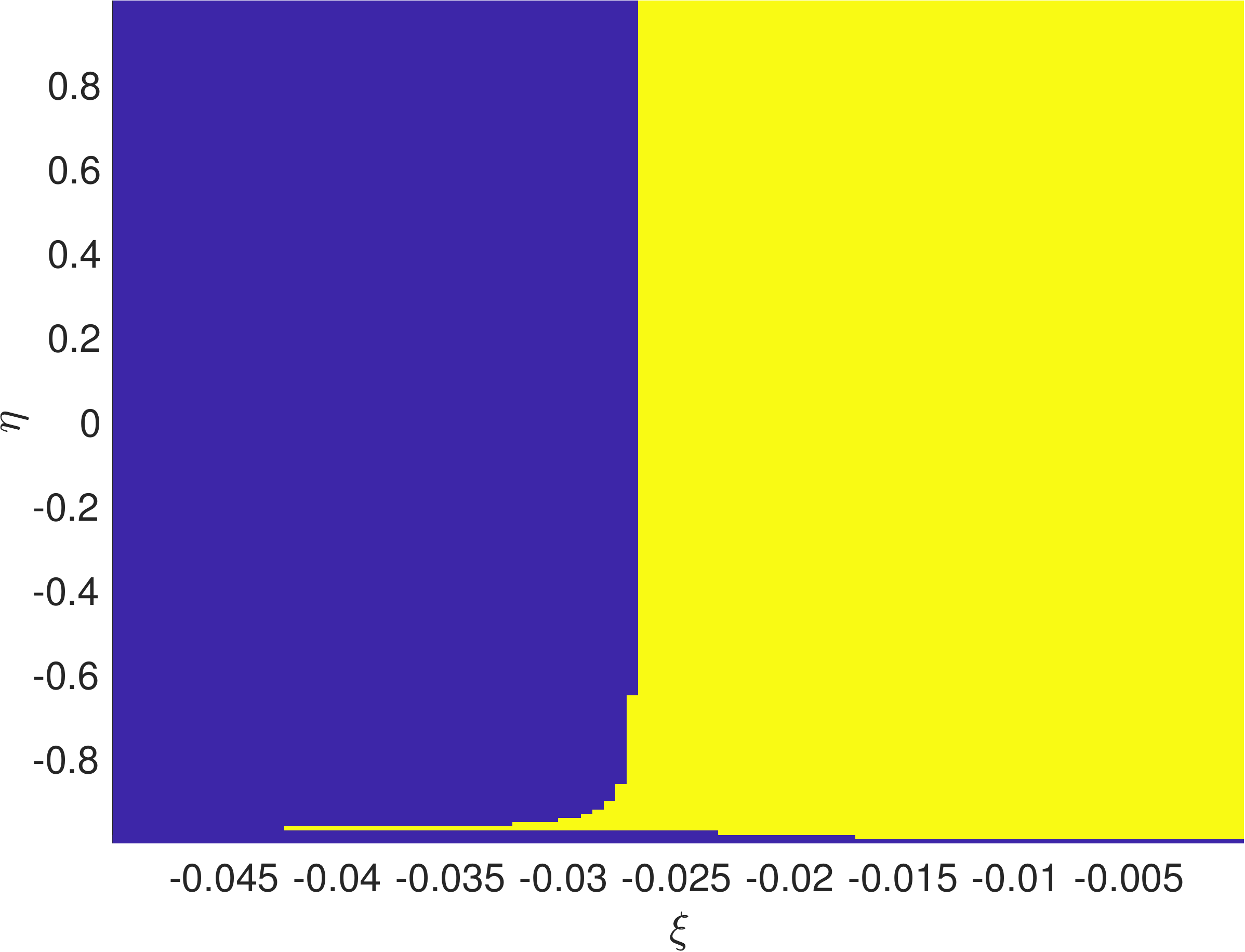}
    \hfill
    \includegraphics [width=.48 \textwidth]{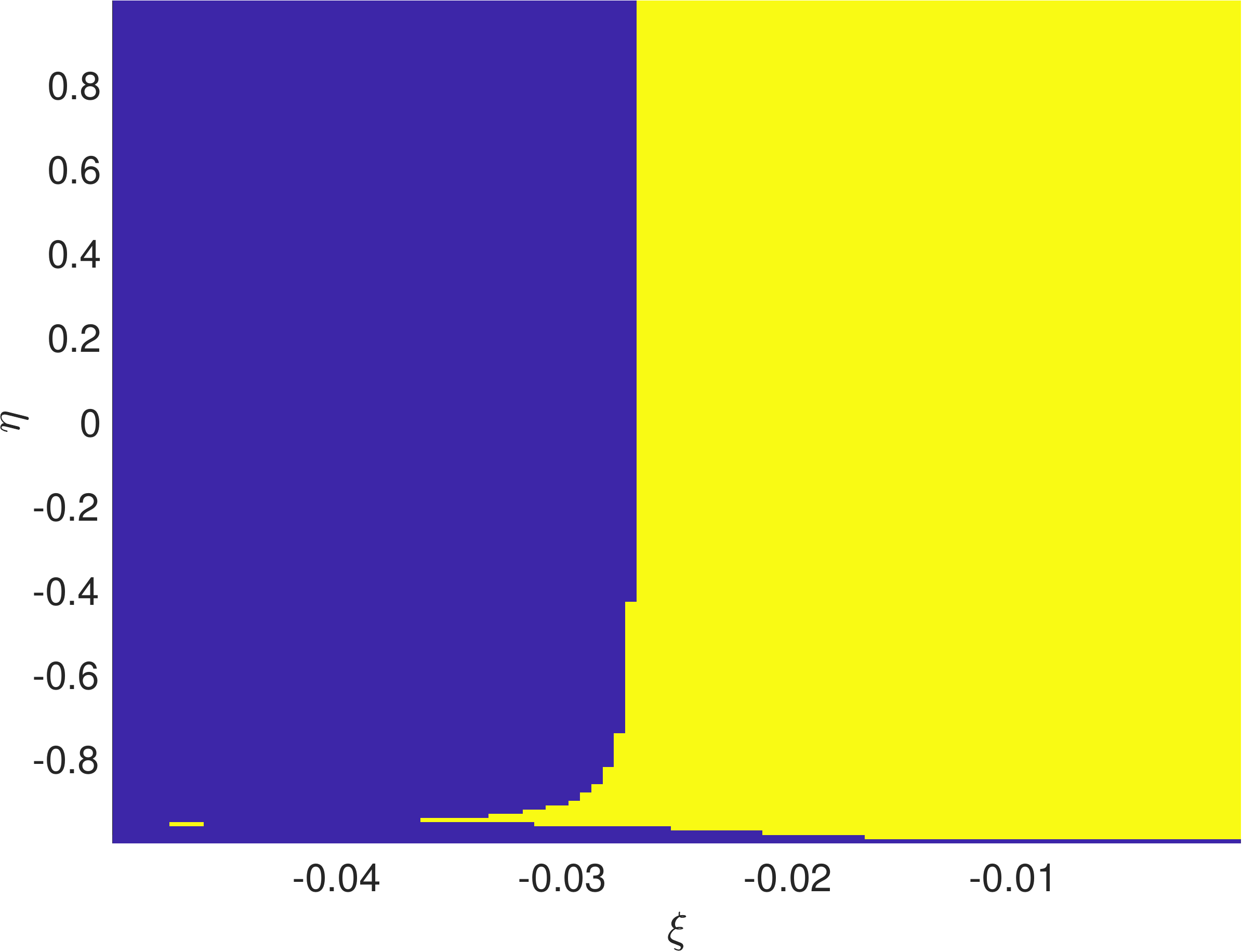}
  }
  \caption{Linear stability plots for the RMIS (left) and MIS
    (right) methods based on the alternate possibility for $T_O$, at
    time-scale separation factors $\kappa=m=10$ (top) and
    $\kappa=m=100$ (bottom), showing the stable region in yellow.  Note that the $\kappa=100$ plots are
    zoomed in to a smaller $\xi$ region than the $\kappa=10$ plots.}
  \label{fig:match2stab}
\end{figure}

To see whether these stability properties translate to other RMIS and
MIS methods (as opposed to only the ones based on the 3/8-Rule), we
additionally examine the linear stability of the RMIS and
MIS methods that use the other 4-stage base method $T_0$ satisfying
the 3rd-order MIS and 4th-order RMIS conditions, given by the
intersection point \eqref{eq:cvals2} from Figure
\ref{fig:RK4stageSol}.  These plots are shown in Figure
\ref{fig:match2stab}.  While the precise shapes of these stability
regions have changed slightly from those in Figure \ref{fig:case1stab}, the
properties observed above are indeed retained for this alternate base
method.

\begin{figure}[tbph!]
  \centerline{
    \includegraphics [width=.49 \textwidth]{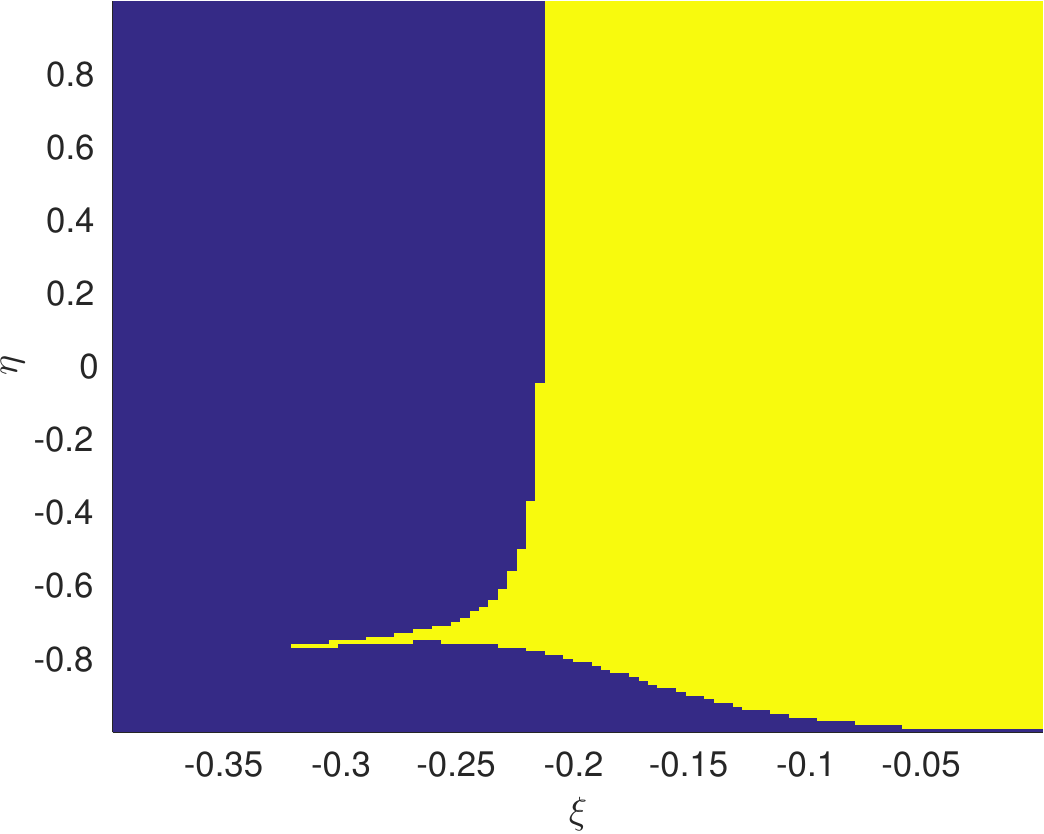}
    \hfill
    \includegraphics [width=.49 \textwidth]{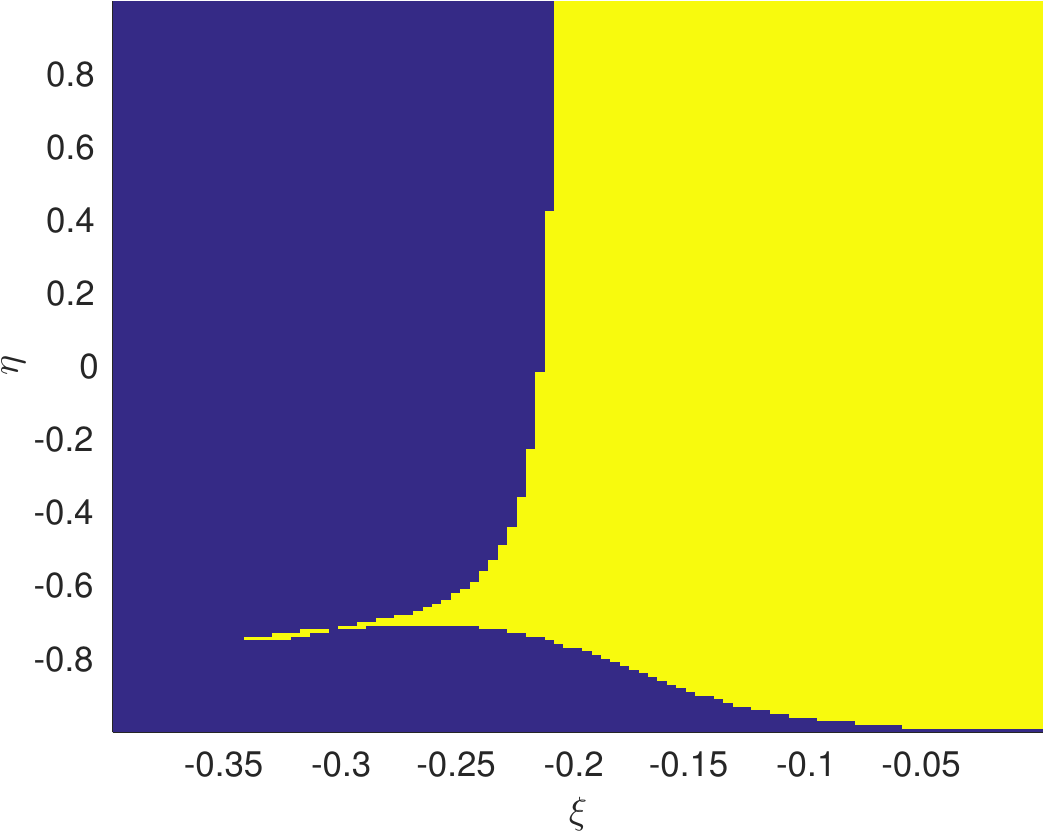}
  }
  \caption{Linear stability plots for the RMIS (left) and MIS
    (right) methods having `maximal' stability regions for time-scale
    separation factor $\kappa=m=10$, showing the stable region in yellow.}
  \label{fig:stabareaOptSingle}
\end{figure}

As a final stability comparison between MIS and RMIS methods, we
investigated the following question: given that there exist
one-parameter families of 3rd-order MIS methods and 4th-order RMIS
methods (as seen in Figure \ref{fig:RK4stageSol}), can we find
the MIS and RMIS methods with `largest' linear stability regions, and
how would these compare against one another?  To this end, we tested
100 base methods along this one-parameter family for each method, and
selected the one with largest stability region (maximal area of the
yellow region).  The results of this investigation for $\kappa=m=10$
are shown in Figure \ref{fig:stabareaOptSingle}.  We again note the
similar stability regions as before, showing slight variations in the
stability boundaries but the same overall shapes.

Based on the results shown in this section, we conclude that the
proposed RMIS methods suffer no deterioration in stability as compared
with their MIS `cousins'.

\section{Convergence and Efficiency Tests}
\label{sec:numerical_tests}

In this section we provide numerical results to verify the analytical
order of the proposed RMIS methods and compare their efficiency
against the MIS methods of \cite{Schlegel2012}.  To that end,
we consider three standard multirate tests, corresponding to
the inverter-chain problem
\cite{Kvaerno1999,Bartel2002,Gunther2001,Constantinescu2013},
a linear multirate problem with strong fast/slow coupling from
Kuhn and Lang \cite{Kuhn2014}, and the well-known brusselator
problem. For each problem we either compare against the analytical
solution or a high-accuracy reference solution.  These reference
solutions were generated by using the $12$th order Gauss implicit
Runge-Kutta method, using fixed time-steps which are $4$ times smaller
than the smallest $h$ value tested for our multirate methods.

We measure solution error in each test by computing the
root-mean-square error over all solution components and over all time
steps $t_n=t_0+hn$,
\begin{equation}
  \label{eq:RMS_error}
  \mathtt{RMSerror} = \left(\frac{1}{MN} \sum_{n=1}^M \left\| y_{n,h}-y_{n,ref} \right\|_2^2\right)^{1/2} ,
\end{equation}
where $y_{n,h}$ and $y_{n,ref}$ are our computed and reference
solutions in $\R^n$ at $t_n$, respectively, $\|\cdot\|_2$ is the
standard vector 2-norm, and $M=(t_{F}-t_{0})/h$ is the total
number of overall time steps of size $h$ taken i n each run over
$t\in\left[t_0,t_F\right]$.  For each test, we present both plots of
$\mathtt{RMSerror}$ versus $h$, as well as tables of the observed
numerical order.  These observed orders are computed using a
least-squares fit to these data where outliers are ignored, i.e., the
numerical orders of convergence only include data where
$10^{-9}\leq\mathtt{RMSerror}\leq 1$.

To assess the efficiency of each method, we utilize standard
error-versus-cost plots.  Since these simulations were performed in
MATLAB, where timings can provide rather poor predictions of runtimes
for true HPC applications, in these tests we measure cost by counting
the total number of ODE right-hand side function calls.  For the
subcycled and telescopic RMIS and MIS methods examined here, this
may be easily calculated as
\begin{equation}
  \label{eq:Total_f_calls}
  \mathtt{TotalFunctionCalls} = s^O + s^I \sum_{i=1}^{s^O} n_i,
\end{equation}
where $n_i$ is the number of fast step subcycles required per slow
stage $i$.  In the context of our five methods (Opt-3/8, RMIS-3/8,
RMIS-KW3, MIS-3/8, and MIS-KW3), we may examine this
accounting for a time-scale separation of $m=100$.  Since the methods
based on the 3/8-Rule perform 34 subcycles per outer stage, and the
methods based on KW3 perform 35 subcycles per outer stage, the total
number of right-hand side function calls per outer step are
$4+4(34+34+34)=412$ and $3+3(35+35+35)=318$.
These per-step costs are summarized in Table \ref{tableMethodEfficiency}.
\begin{table}[!htp]
  \begin{center}
    \begin{tabular}{|c|c|c|c|}
      \hline
      Method Efficiency & $s^O$ & $s^f / s^s $ & $s^O\left(1+\sum_{i=1}^{s^O} n_i\right)$\\
      \hline
      Opt-3/8 & 4 & 102& 412\\
      RMIS-3/8 & 4 & 102 & 412\\
      RMIS-KW3 &  3 & 105& 318\\
      MIS-3/8 & 4 & 102 & 412\\
      MIS-KW3 & 3 &105& 318\\
      \hline
    \end{tabular}
    \caption{Per-step method costs: the choice of $T_O$ determines the
      number of stages $(s^O)$ and subcycles per outer stage $(n_i)$,
      resulting in slight differences in the ratio of fast stages to
      slow stages in each method.
      \label{tableMethodEfficiency}
    }
  \end{center}
\end{table}
With these data, for each test we plot $\mathtt{RMSerror}$ versus
$\mathtt{TotalFunctionCalls}$; ``efficient'' methods correspond to
curves that are nearest the bottom-left corner.

\subsection{Inverter-chain}
\label{ssec:invertDef}

The inverter-chain problem is a partitioned multirate ODE system that
models a chain of MOSFET inverters, which has been used for testing
multirate ODE solvers throughout the literature
\cite{Schlegel2012,Bartel2002,Gunther2001,Constantinescu2013,Savcenco2007,Hundsdorfer2009,Striebel2007,Verhoeven2008,Savcenco2009,Savcenco2010,Oliveira2012}.
The form of the model that we examine is primarily based on the
version used by Kv{\ae}rn{\o} and Rentrop \cite{Kvaerno1999}, although
we utilize an additional scaling term as used in
\cite{Bartel2002,Constantinescu2013}.  The mathematical model is given
by the system of ODEs for $y(t)\in\R^{n_I}$, $0\le t\le 7$:
\begin{equation}
  \label{eq:inverter_chain_model}
  y_k'(t) = y_{op} - y_k(t) - \gamma g_k(t,y), \quad y_k(0)=0, \quad k=1,\ldots,n_I,
\end{equation}
where
\begin{align*}
  g_k(t,y) &= \begin{cases}
    g(y_{in}(t), y_1(t), y_0), & k=1\\
    g(y_{k-1}(t), y_k(t), y_0), & 1<k\le n_I
  \end{cases},\\
  g(y_G, y_D, y_S) &= \left(\max\left(y_G-y_S-y_T,0\right)\right)^2 - \left(\max\left(y_G-y_D-y_T,0\right)\right)^2,
\end{align*}
where $y_{op} = 5$ V, $y_T = 1$ V, $y_0 = 0$ V, $\gamma$ is the
scaling term for tuning the time-scale separation of the problem (we
use $\gamma=100$), and $y_{in}(t)$ is the forcing function,
\[
  y_{in}(t)=\begin{cases}
    0, & 0\le t<5,\\
    t-5, & 5\le t\le 7.
  \end{cases}
\]
We note that $y_{in}(t)$ causes the problem
\eqref{eq:inverter_chain_model} to be non-autonomous.  This is easily
handled since the MIS and RMIS methods are internally consistent, so
we may identify `stage times' $t_{n,l}^{\{q\}} = t_n + c^{\{q\}}_l h$
that correspond to each stage $k_{l}^{\{q\}}$.

\begin{figure}[tbph]
  \centerline{
    \includegraphics[width=.9\linewidth]{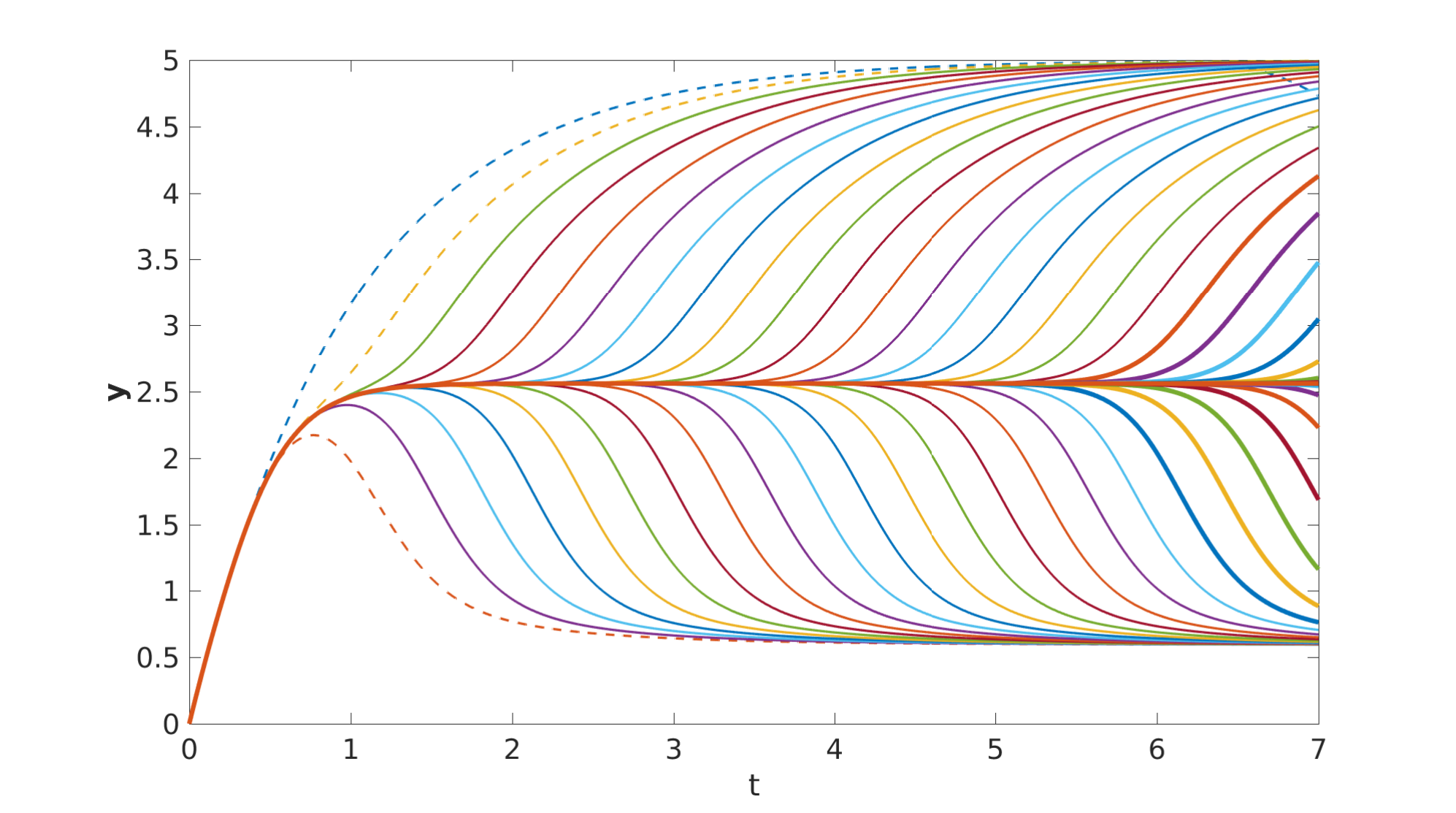}
  }
  \centerline{
    \includegraphics[width=.9\linewidth]{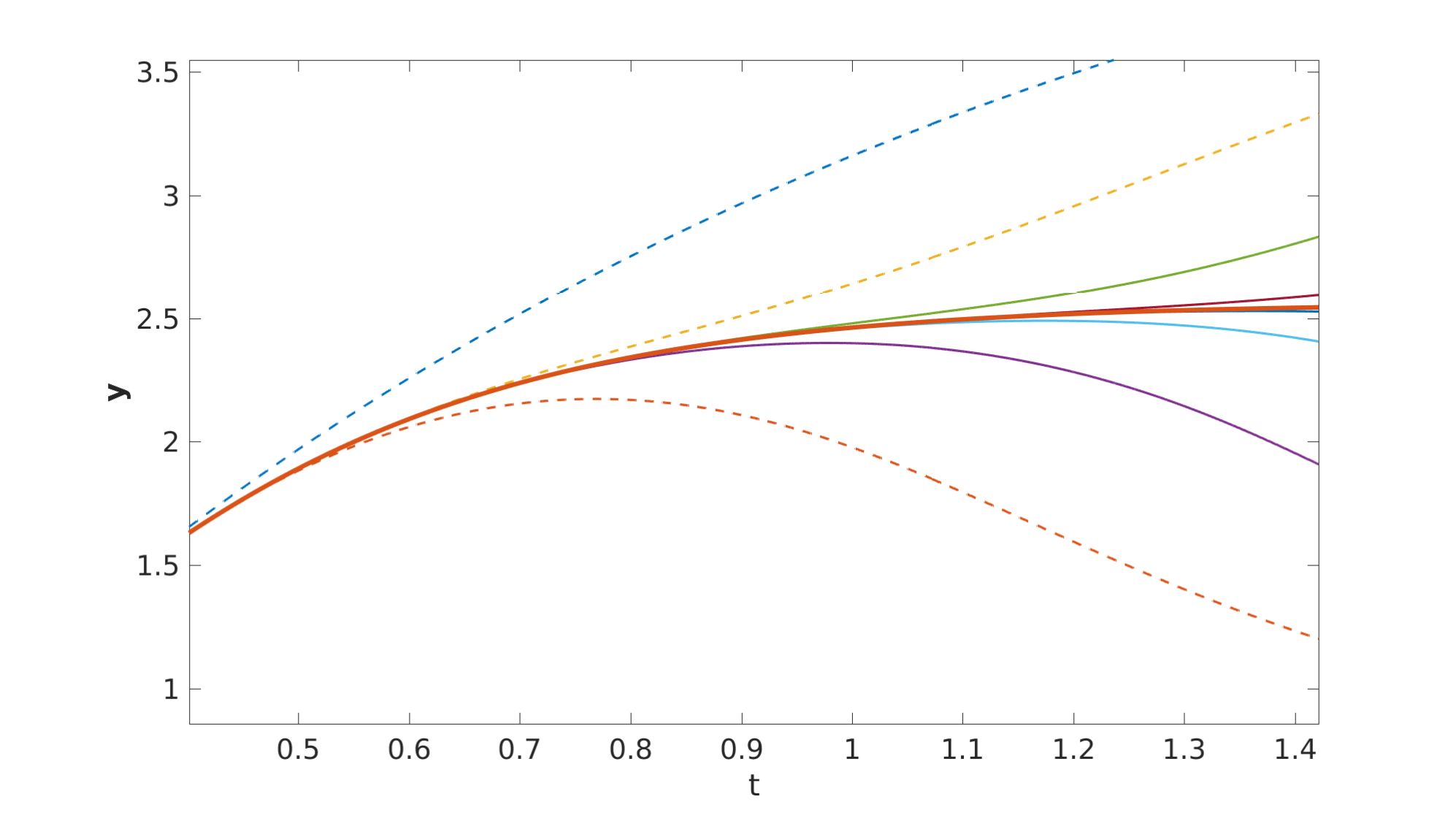}
  }
  \caption{Solutions for the inverter chain problem with $n_I=100$:
    Dotted lines represent the fast components, i.e.~$y_1$ is the blue
    dotted curve at the top-left and $y_{100}$ remains at the value
    2.5 at the final time.  At bottom is a zoom-in of the initial
    departure of fast components, indicating that even small
    differences in the initial integration solution may result in
    disparate overall solution values.}
  \label{fig:inverterSol}
\end{figure}

The time scales for these inverters decrease with index, so we
partition this so that the first $b$ equations are ``fast'' and the
remainder are ``slow'',
\begin{align*}
  \flam{f}(t,y)^{\intercal} &= \begin{bmatrix} y_1'(t) & \cdots& y_b'(t)& 0& \cdots& 0\end{bmatrix}\\
  \flam{s}(t,y)^{\intercal} &= \begin{bmatrix} 0& \cdots& 0& y_{b+1}'(t)& \cdots& y_{n_I}'(t)\end{bmatrix}.
\end{align*}
In our tests, we use $n_I=100$ total inverters, with the first $b=3$
grouped into $\flam{f}(t,y)$.  We note that in other papers that use a
similar setup, the first $b=20$ were chosen as ``fast''; however in
our tests we found that only 3 were required for our desired multirate
time-scale separation factor of $m=100$.
In Figure \ref{fig:inverterSol} we show
the solutions for this problem, as well as a zoom-in of the initial
departure of the fast inverters from the larger group.  In our tests
we found that accuracy in these initial departures proved crucial for
overall solution accuracy.

\begin{figure}[tbph]
  \centerline{
    \includegraphics[width=.48\linewidth]{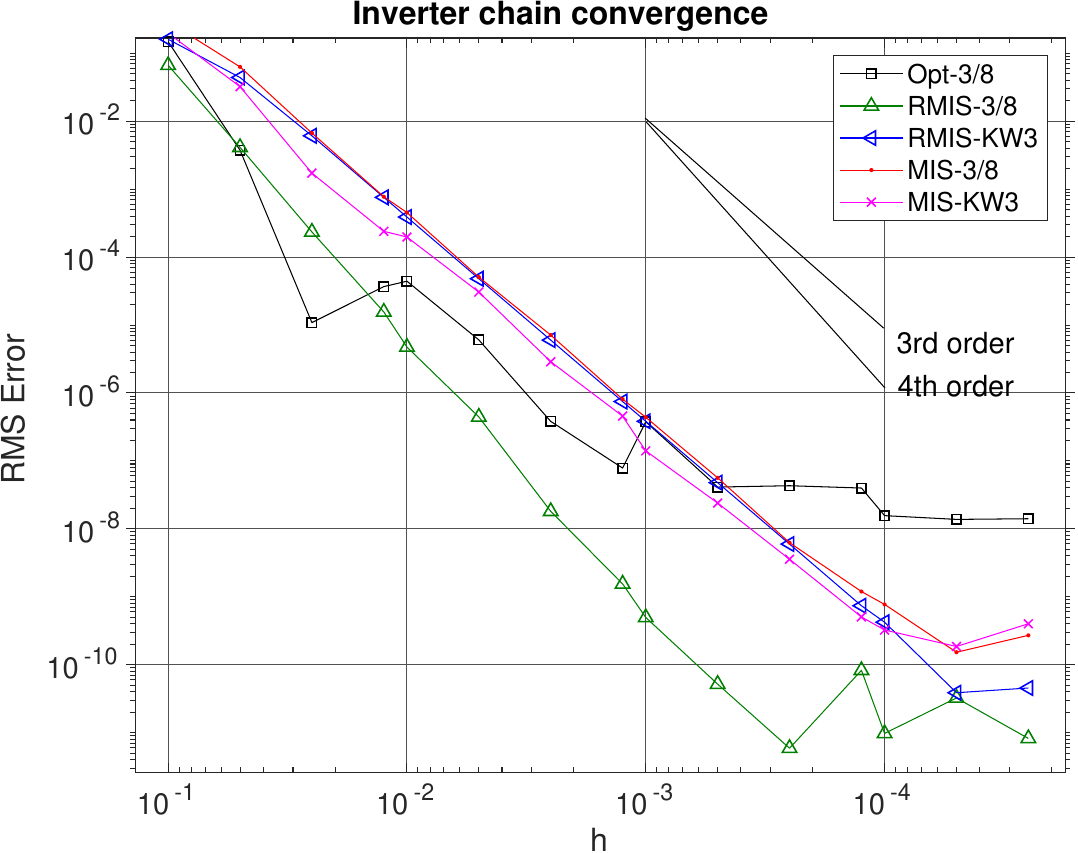}
    \hfill
    \includegraphics[width=.48\linewidth]{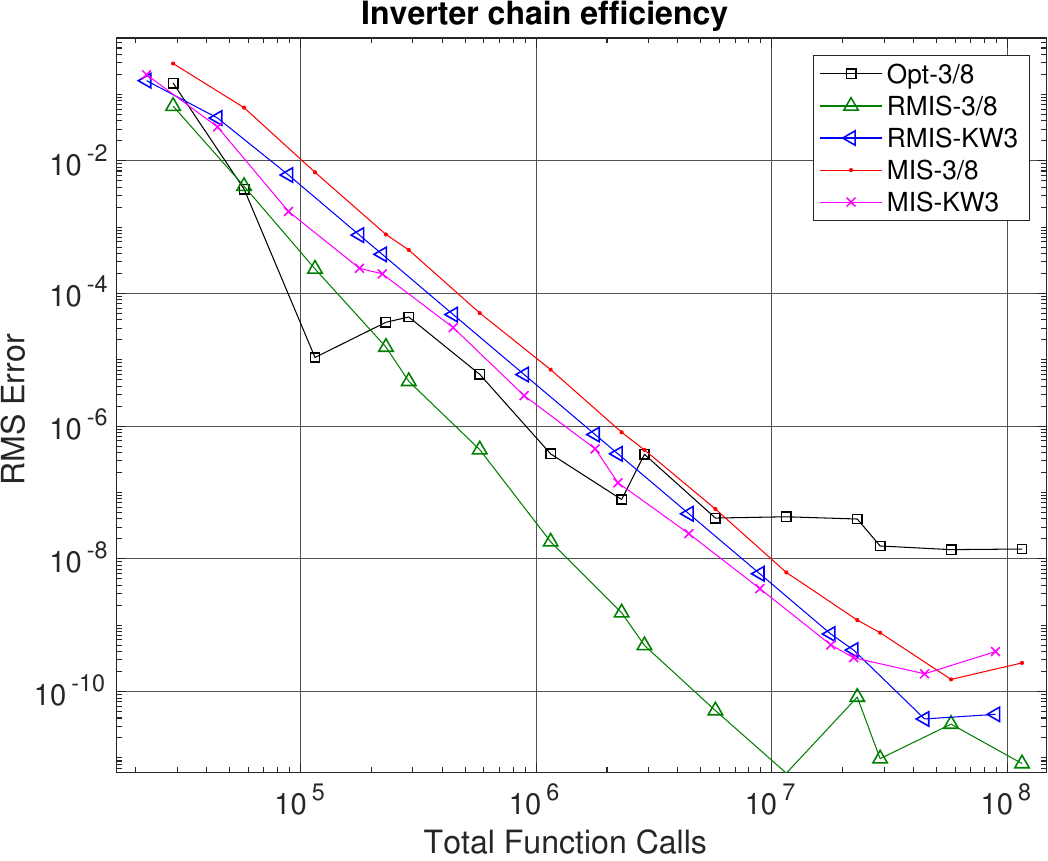}
  }
  \caption{Convergence (left) and efficiency (right) for the inverter
    chain problem: the convergence results are consistent with
    expectations for all methods except Opt-3/8, where the early
    flattening of the error is likely due to increased sensitivity to
    floating-point round-off.  Here, the most efficient methods for
    larger error values are RMIS-3/8 and Opt-3/8, while for smaller
    errors RMIS-3/8 is the clear winner; all other methods perform
    comparably well.}
  \label{fig:inverterResults}
\end{figure}

In the left portion of Figure \ref{fig:inverterResults} we plot
$\mathtt{RMSerror}$ versus step size $h$ for each of the five methods
tested.  We highlight a few observations in these results.  First,
all methods demonstrate convergence as $h\to 0$ to a point, beyond
which convergence stagnates.  This stagnation point is below $10^{-9}$
for all but the Opt-3/8 method, that stagnates slightly earlier at
around  $10^{-8}$, indicating a reference solution accurate to
approximately $10^{-9}$ for this test.  We hypothesize that the
Opt-3/8 method is more susceptible to accumulation of floating-point
round-off than the other methods, since unlike the others that are
defined by a structured pattern of small-valued coefficients, Opt-3/8
has more widely-varying coefficients $\lam{b}{f}_i \in
\left[-5.5\times 10^4, 1.2\times 10^5\right]$.
Second, both the RMIS-3/8 and Opt-3/8 methods show faster rates of
convergence than the other methods; unfortunately, the stagnation of
Opt-3/8 halts this fast convergence somewhat early, but RMIS-3/8 shows
consistently fast convergence until below the reference solution
accuracy.  The best-fit orders of convergence for the results from
this figure are: 1.74 (Opt-3/8), 4.07 (RMIS-3/8), 2.93 (RMIS-KW3),
2.98 (MIS-3/8) and 2.98 (MIS-KW3). We note that all
show their expected rates of convergence except for Opt-3/8, which is
likely due to its larger error floor, since prior to that point it is
converging at least as rapidly as the 3rd-order methods.

Accuracy alone provides only an incomplete picture of performance,
since the methods using the 3/8-Rule require 25\% more function calls
per step than those using KW3.
To this end, in the right portion of
Figure \ref{fig:inverterResults} we also plot the $\mathtt{RMSerror}$ versus
$\mathtt{TotalFunctionCalls}$ for each of the five methods.  Although
the blue and magenta curves for the KW3-based methods indeed shift
further to the left in relation to the other methods, RMIS-3/8 is
still the most efficient of all the methods at nearly all error
values, and is only outperformed by Opt-3/8 for relatively loose error
values (above $\sim\! 10^{-4}$).  The efficiency benefit of RMIS-3/8
is most notable for errors below $\sim\! 10^{-8}$, where it requires
approximately 10 times less work than the other methods.


\subsection{Strongly-coupled Linear Test}
\label{ssec:kuhnDef}

\begin{figure}[tbph]
  \centerline{
    \includegraphics[width=.95\linewidth]{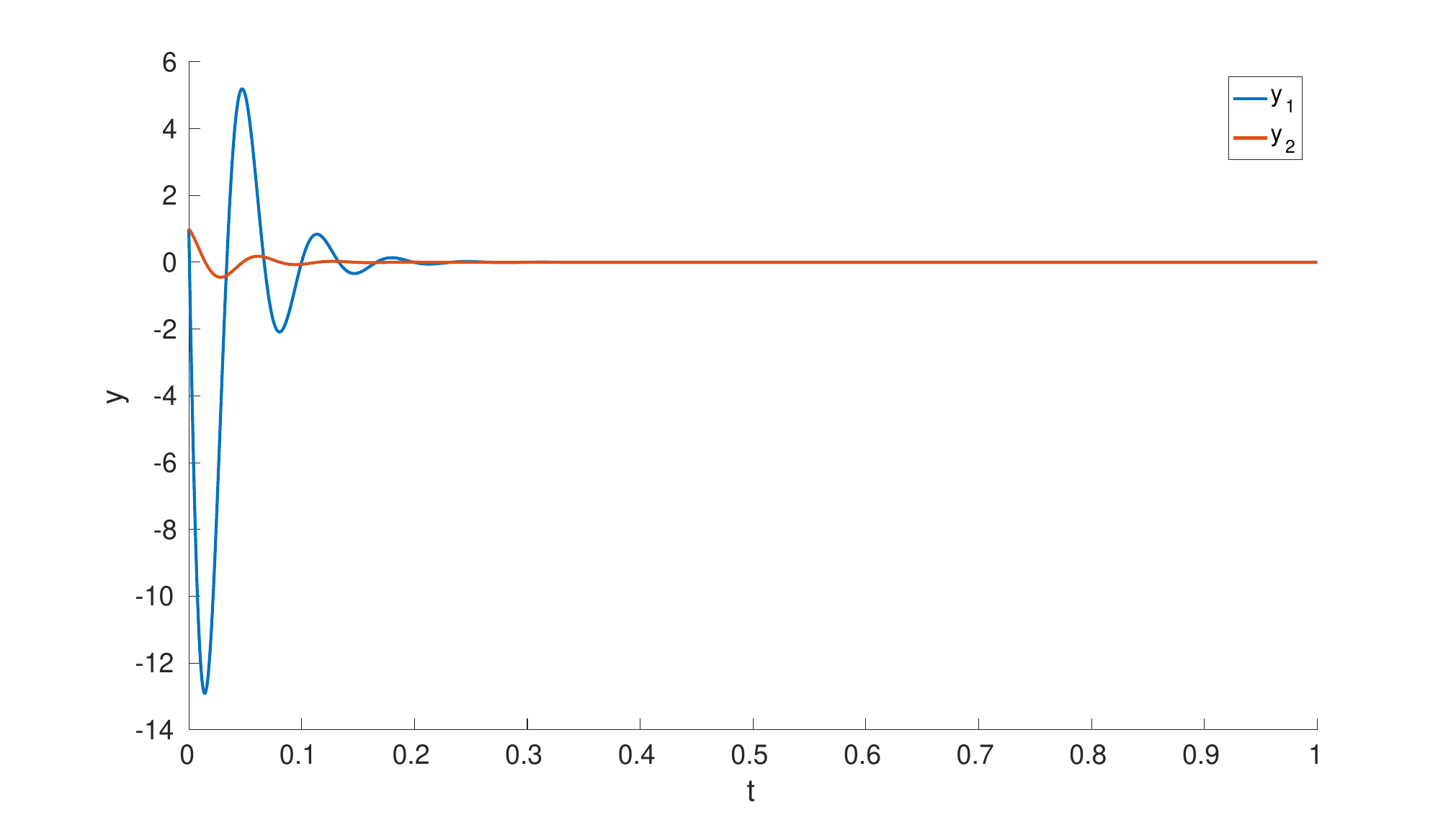}
  }
  \caption{Solutions for the linear test problem: note that the
    ``slow'' component $y_2$ varies far less than strongly than the
    ``fast'' component $y_1$.}
  \label{fig:kuhnSol}
\end{figure}

As a second test, we use a linear ODE system with strong fast/slow
coupling that was used by Kuhn and Lang in their studies of multirate
stability \cite{Kuhn2014}; variants of this problem have been used by
a variety of authors in testing multirate algorithms
\cite{Kvaerno2000,Constantinescu2013,Kuhn2014}.  Our motivation for
this test is to more rigorously explore accuracy and efficiency of the
RMIS and MIS algorithms in the face of strongly-coupled problems,
thereby exercising the coupling matrices $\lam{A}{f,s}$ and
$\lam{A}{s,f}$, and exploring whether the extreme sparsity of our
$\lam{b}{f}$ in the RMIS algorithm causes trouble for strongly-coupled
problems.  The ODE system is identical to the partitioned test
problem \eqref{eq:partitioned_linear} examined in Section
\ref{sec:linear_stability}, where here the system-defining matrix is
given by
\[
  G = \begin{bmatrix} -5 & -1900\\ 5 & -50\end{bmatrix}
\]
and the problem is evolved over the time interval $t\in [0,1]$.  The
eigenvalues of $G$ are complex conjugates,
\[
  \lambda = - \frac{55\pm 5i\sqrt{1439}}{2},
\]
giving rise to the analytical solution
\[
  \begin{bmatrix}  y_1\\  y_2 \end{bmatrix}
  =
  e^{-\frac{55}{2}t}
  \begin{bmatrix}
    \cos\left(\frac{5\sqrt{1439}t}{2}\right)-\frac{751}{\sqrt{1439}}\sin\left(\frac{5\sqrt{1439}t}{2}\right)\\
    \cos\left(\frac{5\sqrt{1439}t}{2}\right)-\frac{7}{\sqrt{1439}}\sin\left(\frac{5\sqrt{1439}t}{2}\right)
  \end{bmatrix}.
\]
which is shown in Figure \ref{fig:kuhnSol}.  The time-scale separation
in this problem arises from the strong fast/slow coupling, rendering
the `sin' term in $y_1$ approximately 100 times stronger than in
$y_2$, at least until the solutions decay to zero.

We cast this problem into additive multirate form by
partitioning the right-hand side into fast and slow components,
\[
  \flam{f}(y) = \begin{bmatrix} -5 & -1900\\ 0 & 0\end{bmatrix} \begin{bmatrix} y_1\\y_2 \end{bmatrix}
  \quad\text{and}\quad
  \flam{s}(y) = \begin{bmatrix} 0 & 0\\ 5 & -50\end{bmatrix} \begin{bmatrix} y_1\\y_2 \end{bmatrix}.
\]
\begin{figure}[tbph]
  \centerline{
    \includegraphics[width=.48\linewidth]{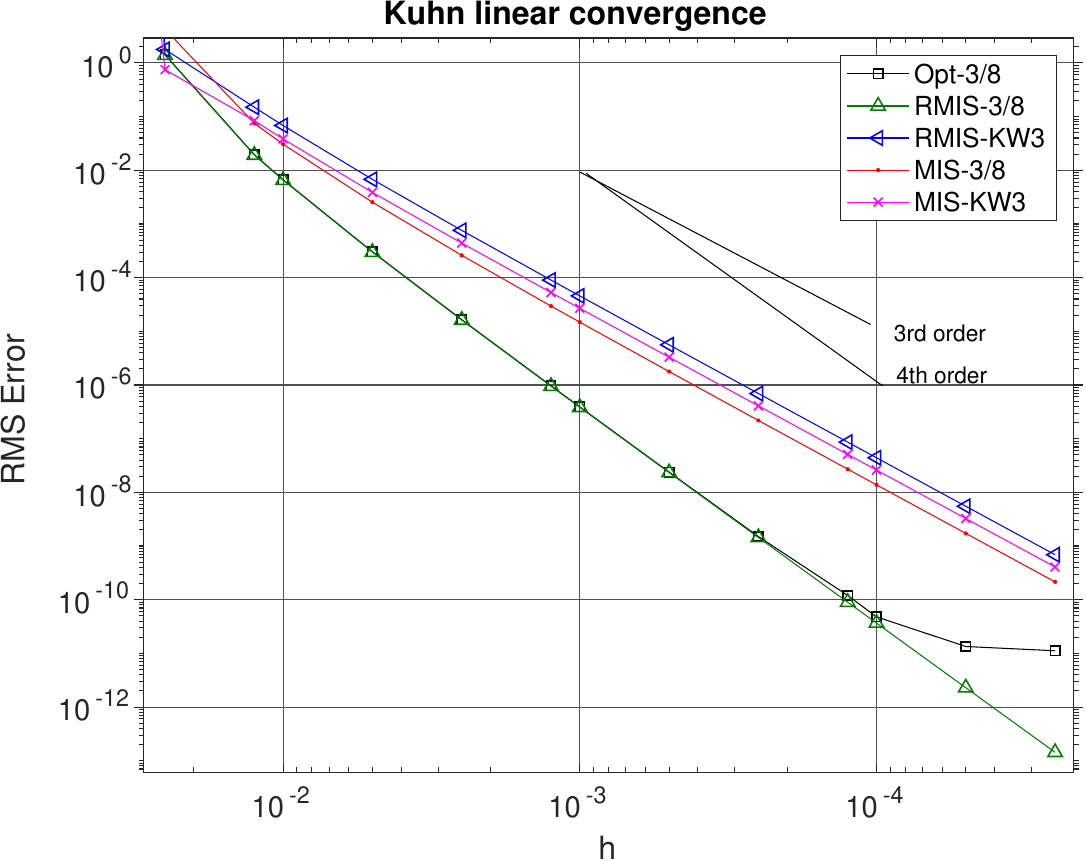}
    \hfill
     \includegraphics[width=.48\linewidth]{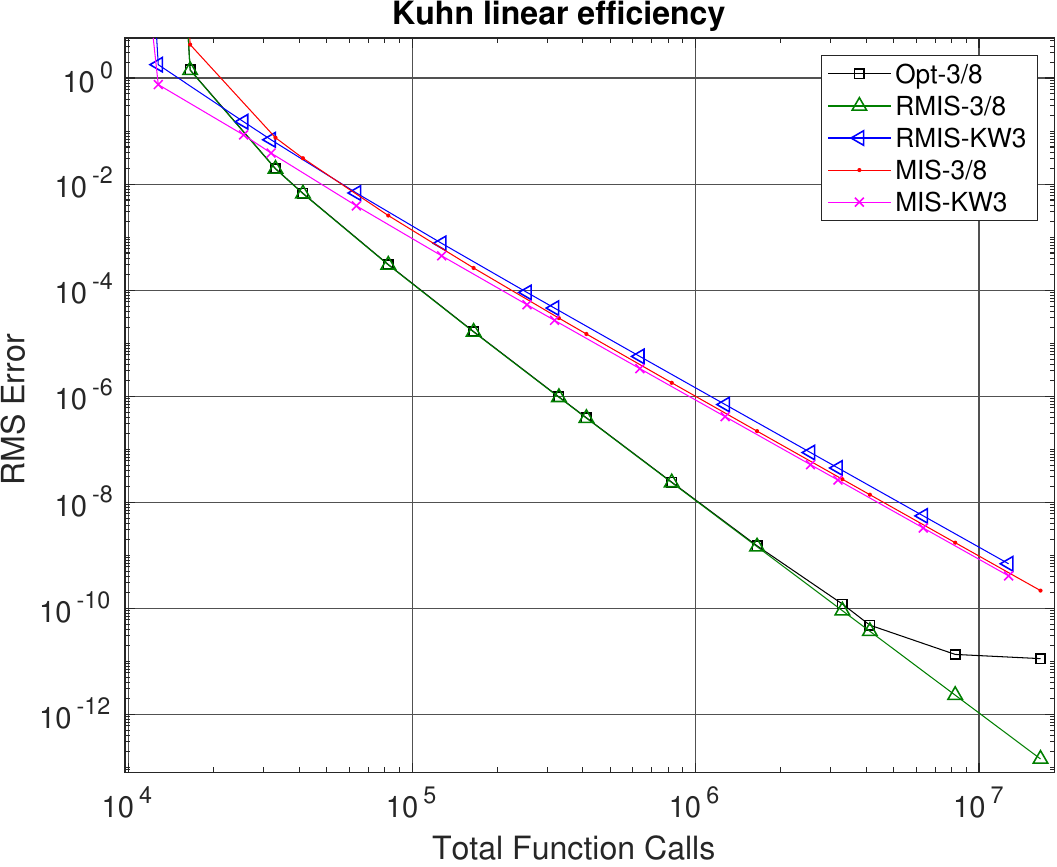}
  }
  \caption{Convergence (left) and efficiency (right) for the linear
    test problem: the convergence results are consistent with
    expectations for all methods, and the 4th-order methods show
    significantly better efficiency than the 3rd-order methods.}
  \label{fig:kuhnResults}
\end{figure}

In Figure \ref{fig:kuhnResults} we plot both $\mathtt{RMSerror}$
versus $h$ (left) and $\mathtt{RMSerror}$ versus
$\mathtt{TotalFunctionCalls}$ (right).  This problem exhibits very
straight order of accuracy lines for all methods except Opt-3/8, which
begins to level off at approximately $10^{-11}$, again likely due to
increased sensitivity to accumulation of floating-point round-off
errors.  The best-fit convergence orders for these results are: 4.22
(Opt-3/8), 4.22 (RMIS-3/8), 3.09 (RMIS-KW3), 3.18 (MIS-3/8) and 3.09
(MIS-KW3).  We do not fully understand the higher-than-expected
orders of accuracy for Opt-3/8 and RMIS-3/8, except that since the
problem is linear then many of the GARK order conditions are
irrelevant.  As with the previous problem, the reduced cost per
step of the KW3-based methods is insufficient to surpass the
efficiency of the 4th-order methods, with Opt-3/8 and RMIS-3/8
demonstrating essentially-identical efficiency for errors larger than
$10^{-10}$, and with RMIS-3/8 proving more efficient past that point.



\subsection{Brusselator}
\label{ssec:brussDef}

\begin{figure}[tbph]
  \centerline{
    \includegraphics[width=.9\linewidth]{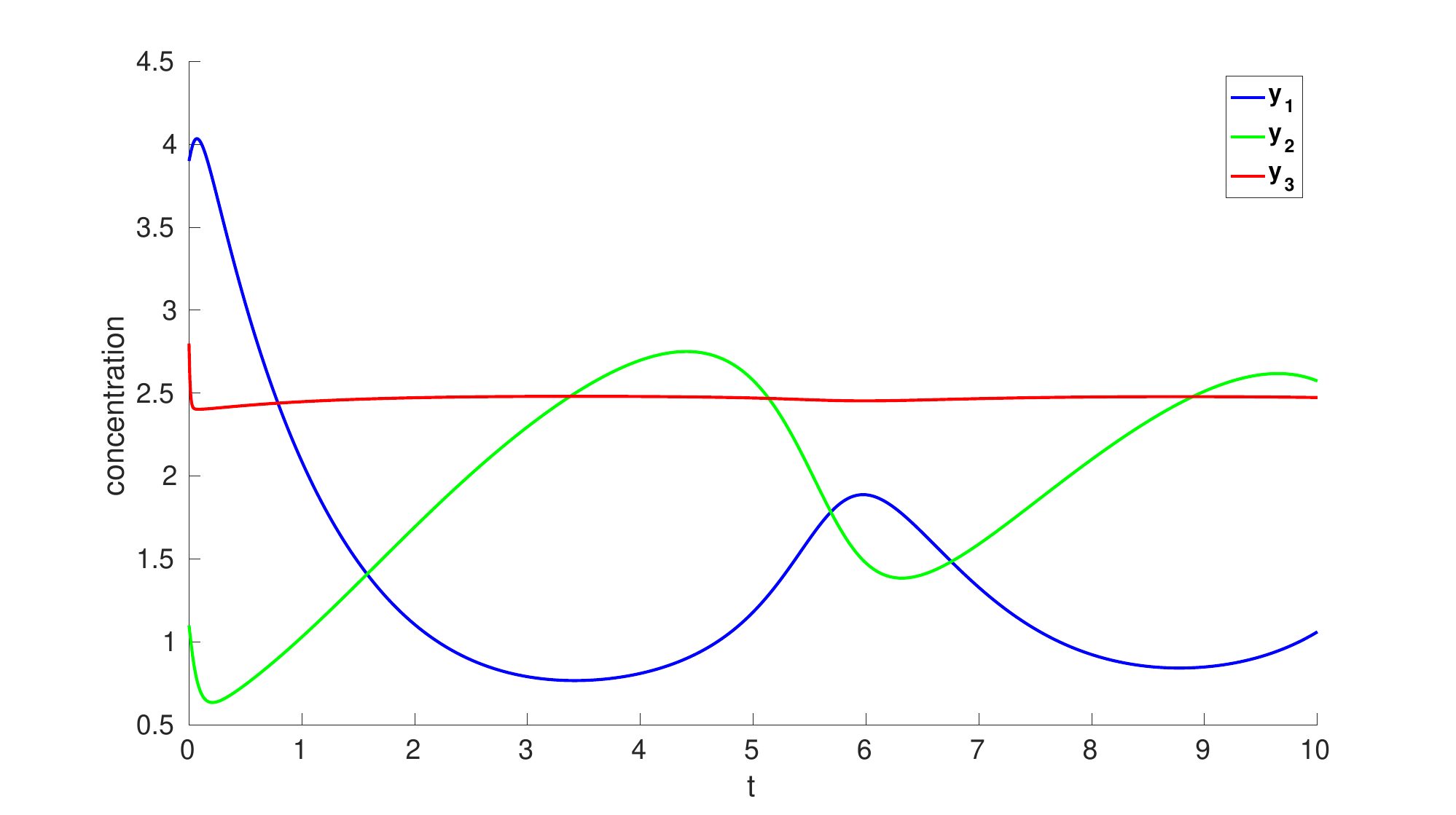}
  }
  \caption{Solutions for the Brusselator test problem
    \eqref{eq:brusselator}: note that all components vary rapidly at
    first and only slowly thereafter.}
  \label{fig:brussSol}
\end{figure}

We consider a system of stiff nonlinear ODEs that captures some of the
physical challenges of the brusselator chemical reaction network
problem, first described as a 1D PDE by Prigogine in 1967
\cite{Prigogine1967}.  More recently this test was used in a
multiphysics paper by Estep in 2008, and a general
computational multiphysics review paper in 2013
\cite{Estep2008,Keyes2013}.  Our version of the problem is a tunable
two-rate initial-value problem represented as a system of three
nonlinearly-coupled ODEs,
\begin{align}
  \notag
  y'(t) &= \flam{f}(y)+\flam{s}(y), \quad 0\le t\le 10\\
  \notag
  y(0) &= \begin{bmatrix} 3.9 & 1.1 & 2.8 \end{bmatrix}^{\intercal}\\
  \label{eq:brusselator}
  \text{where}\qquad &\\
  \notag
  \flam{f}(y) &= \begin{bmatrix} 0\\ 0\\ \frac{b-y_3}{\varepsilon} \end{bmatrix}
  \quad\text{and}\quad
  \flam{s}(y) = \begin{bmatrix} a - (y_3+1)y_1 + y_2y_1^2\\ y_3y_1 - y_2y_1^2\\ -y_3y_1 \end{bmatrix},
\end{align}
where the parameters are chosen to be $a=1.2$, $b=2.5$, and
$\varepsilon=10^{-2}$.  As shown in the above partitioning, the fast
function $\flam{f}$ contains only the term which is scaled by
$\varepsilon$; the problem time-scale separation is approximately
$1/\varepsilon=100$.  With this particular setup, the problem exhibits
a rapid change in the solution at the start of the simulation for
$t<0.2$, with slower variation for the remainder of the time interval,
as shown in Figure \ref{fig:brussSol}.

\begin{figure}[tbph]
  \centerline{
    \includegraphics[width=.48\linewidth]{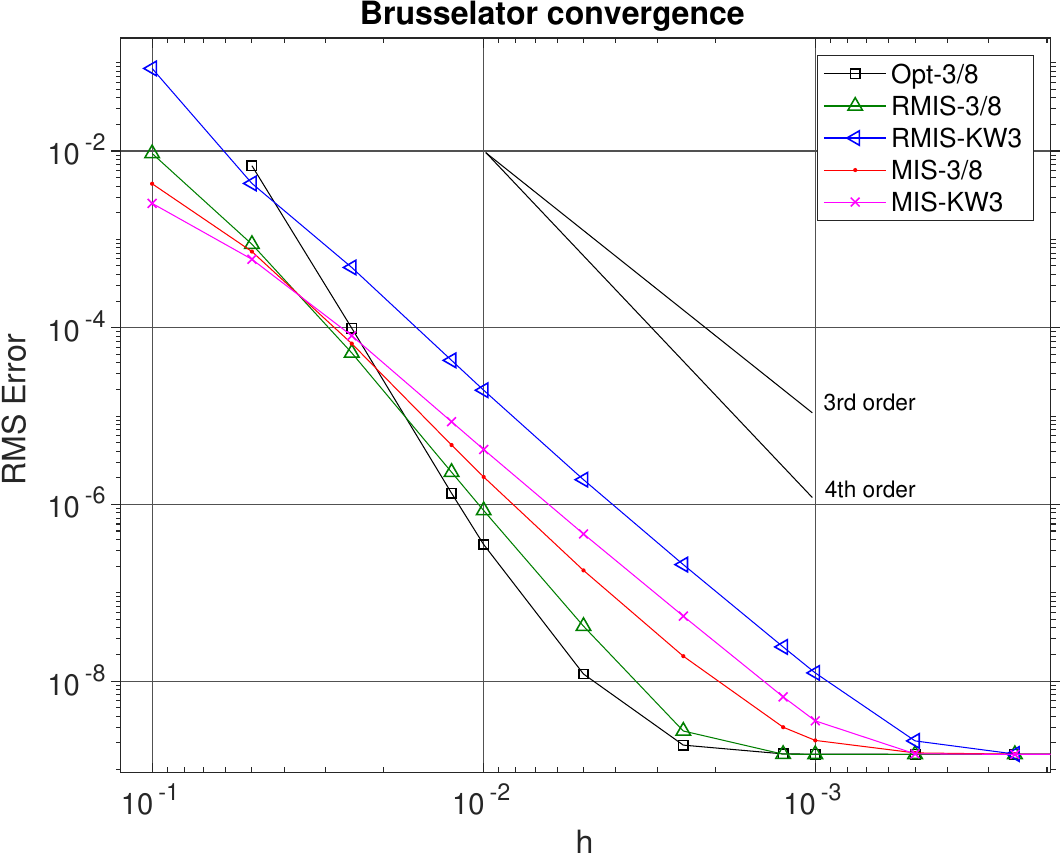}
    \hfill
    \includegraphics[width=.48\linewidth]{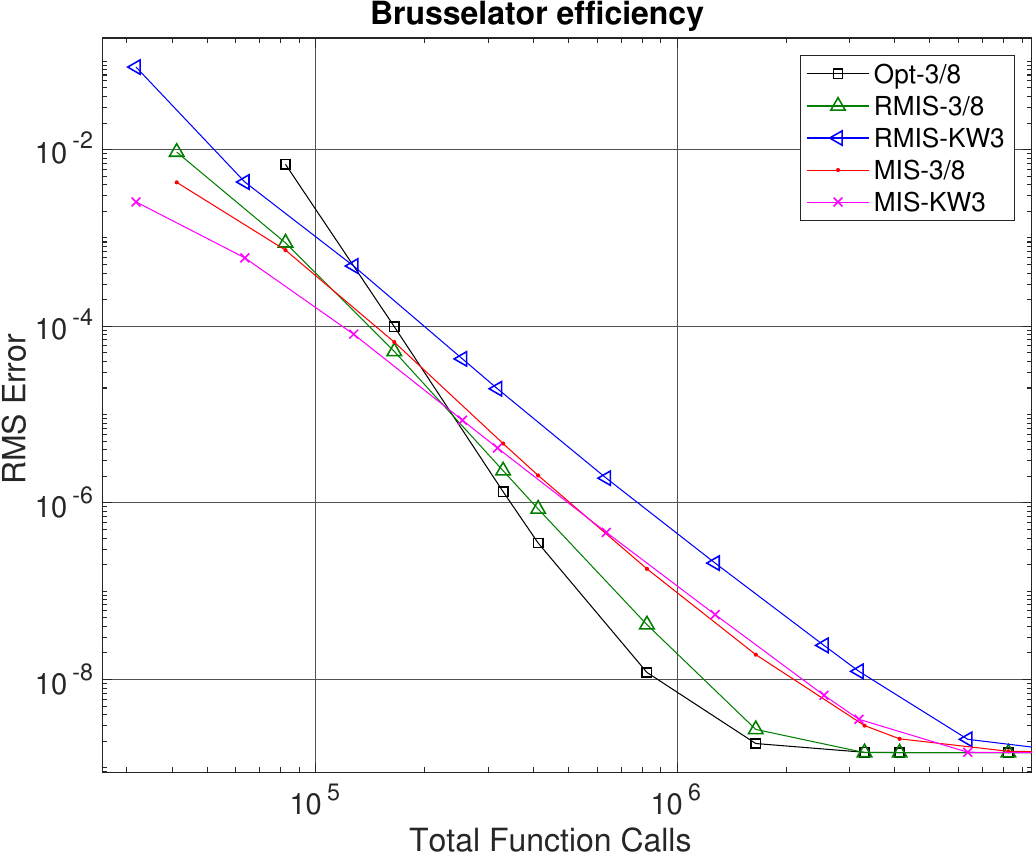}
  }
  \caption{The optimized method and the RMIS seem to do best. Method
    efficiency for this problem is dependent on both the order of
    accuracy and on the constant multiplier on the error terms}
  \label{fig:brusselatorResults}
\end{figure}

In the left portion of Figure \ref{fig:brusselatorResults} we plot
$\mathtt{RMSerror}$ versus $h$, and on the right we plot
$\mathtt{RMSerror}$ versus $\mathtt{TotalFunctionCalls}$.
These plots show perhaps the most interesting results of the three
tests.  Focusing first on the convergence plots we note an
intersection point in the curves, showing that the optimal choice of
method can depend more intimately on the desired solution error.
Beginning at the large-error end, we see that the MIS-KW3 method has
smallest error at the largest tested step sizes, although those errors
are not significantly different than those for MIS-3/8 or RMIS-3/8.
At smaller $h$ and smaller error values, the increased accuracy of
RMIS-3/8 and Opt-3/8 appear, rapidly achieving errors $\sim\! 100$ times
smaller than the 3rd-order methods at step sizes $\sim\! 2\times
10^{-3}$.  Convergence of all methods stagnates at approximately
$10^{-9}$, clearly indicating the accuracy of our reference solution.
The corresponding best-fit estimates of the convergence orders for
these results are: 5.83 (Opt-3/8), 4.16 (RMIS-3/8), 3.30 (RMIS-KW3),
3.28 (MIS-3/8) and 3.02 (MIS-KW3).  We point out the
`superconvergent' behavior of Opt-3/8, which suggests that the
optimization successfully minimized the dominant 5th-order error terms
of relevance to this problem. The other methods have an observed
numerical order of accuracy consistent with our theoretical
expectations.

When examining the efficiency plot at the right of Figure
\ref{fig:brusselatorResults}, we observe that the decreased cost per
step of the KW3-based methods shift those results to the left; making
MIS-KW3 the clear winner for larger error values.  However, at error
values below $\sim\! 10^{-5}$, the higher order methods begin to
outperform the lower order methods, with Opt-3/8 the most efficient,
followed by RMIS-3/8.

\section{Conclusions}
\label{ch:Conclusions}

In this work, we propose a variation of the existing MIS
multirate methods \cite{Schlegel2012, Schlegel2012a, Schlegel2012b,
  Knoth2012, Knoth2014, Knoth1998, Schlegel2009, Wensch2009}.
Both the original MIS methods, as well as our extensions, demonstrate
a number of very attractive properties for multirate integration.
These methods are telescopic, allowing for recursion to any number of
problem time scales, allow subcycling of the fast method within the
slow, and are highly flexible in that they allow for varying the
time-scale separation of the method ($m$) between steps.  Finally,
when constructed using inner and outer Runge--Kutta methods, $T_I$ and
$T_O$, of at least third order, MIS methods are at worst second order
accurate and achieve third order when $T_O$ satisfies the auxiliary
condition \eqref{eq:rfsmr_3rd}.

Recent theory by Sandu and G{\"u}nther complements these MIS methods
nicely.  Specifically, they have proposed a general formulation for
analyzing a wide range of Runge-Kutta-like methods, named
``Generalized-structure Additive Runge-Kutta'' (GARK) methods
\cite{Sandu2015}, that lays a strong theoretical foundation for
understanding and extending MIS methods.  More beneficial, however,
is their subsequent work that directly analyzed MIS methods using
their GARK framework \cite{Gunther2016}, opening the door for
subsequent extensions to this analysis.

In this context, we propose a new multirate algorithm which we name
``relaxed multirate infinitesimal step'' (RMIS), due to a relaxation
from the MIS approach on how the fast stage solutions are
combined to construct the overall multirate step solution (i.e., the
coefficients $\lam{b}{f}$).  This simple change, while leaving the
remainder of the MIS approach intact (i.e., the algorithmic approach
as well as the coefficients $\lam{A}{f,f}$, $\lam{A}{f,s}$,
$\lam{A}{s,f}$, $\lam{A}{s,s}$ and $\lam{b}{s}$), allows for a number
of remarkable extensions to their work.  First, we are able to
construct up to fourth order, and at least third order, multirate
methods with comparable stability and improved efficiency to MIS
methods, where the determining factors for order four are: $T_I$ must
be at least third order and have explicit first stage; $T_O$ must be
explicit, at least fourth order, and must satisfy the condition
\eqref{eq:rmis_4th}.  A key component of this analysis is our result
from Lemma \ref{lem:RMIS_bf_simplifications}, where we show that
due to our selection of the coefficients $\lam{b}{f}$, the fast order
conditions are equivalent to the slow conditions, so one only needs to
prove half of the conditions to guarantee overall order of
accuracy.  These new methods show nearly identical stability
properties, and retain all of the above ``attractive properties''
described above for MIS methods (telescopic, subcycling-ready,
flexible to adjusting $m$).  In addition, since MIS and RMIS only
differ in their selection of the $\lam{b}{f}$ coefficients, then it is
straightforward to define MIS as an embedding within an RMIS method.

In addition to providing detailed proof of the fourth order conditions
for the RMIS algorithm, and comparison of the linear stability between
MIS and RMIS methods, we provide numerical comparisons of the
performance of multiple RMIS and MIS methods on three standard
multirate test problems in Section \ref{sec:numerical_tests}.  These
problems include the standard ``inverter-chain'' system of nonlinearly
coupled ODEs, a linear multirate problem with strong fast/slow
coupling, and the standard ``brusselator'' stiff nonlinear ODE
system.  Through these experiments, we note that our theoretical
expectations regarding the order of accuracy for each method are borne
out in the results for each problem, namely that the RMIS methods
satisfying \eqref{eq:rmis_4th} exhibit fourth order convergence while
those that do not are only third order accurate, and that the MIS
methods are at best third order accurate (even when using higher-order
base methods).  As a result, due to its consistently strong
performance (stability, error and efficiency), the proposed RMIS-3/8
method is a clear contribution to the ever-growing `stable' of
multirate methods, particularly when highly-accurate solutions are
desired.

Although at the time when we performed this work no other
fourth-order multirate methods existed, the subsequent release of
preprints by Sandu and collaborators demonstrating alternate
approaches for fourth-order multirate methods warrants discussion
here.  Each of the available fourth-order methods (MrGARK \cite{Sarshar2018},
MRI-GARK \cite{Sandu2018}, and RMIS), follow distinct approaches for
obtaining higher order.  The MrGARK methods pre-select the tables
$\lam{A}{f,f}$ and $\lam{A}{s,s}$ and the value $m$, and then
construct fast/slow coupling coefficients $\lam{A}{f,s}$ and
$\lam{A}{s,f}$ to satisfy requisite coupling conditions.  The MRI-GARK
and RMIS methods, on the other hand, define a set of modified ``fast''
initial value problems using the slow right-hand side function values,
$\flam{s}$, allowing for significantly more flexibility in the choice
of integrator for the fast time scale.  However, these approaches
differ in two fundamental areas.  First, the fast time-scale IVP
modification in RMIS is simpler than that for MRI-GARK, requiring
storage of only a single vector to encode the slow-to-fast source
terms, as opposed to a vector-valued time-dependent function.  While
perhaps not as costly as other components of the solves, this reduced
memory and computation footprint hints that RMIS may be more efficient
per-step than MRI-GARK.  Second, the approach for obtaining fourth
order differs dramatically between these methods; RMIS methods will
automatically obtain this accuracy through proper selection of the
outer Butcher table, $T_O$ -- as shown in Figure \ref{fig:RK4stageSol}
there are infinitely many options in this regard.  MRI-GARK methods,
on the other hand, achieve higher order through appropriate selection
of the polynomial coefficients $\gamma_{i,j}\left(\tau\right)$ in
equation \eqref{eq:mri_gark_fast_ode}; these coefficients must
themselves be carefully chosen to enforce a modified set of order
conditions. Although fourth-order RMIS methods are thus more flexible
with respect to the fast time scale than MrGARK, and are both more
efficient and easier to construct than MRI-GARK, it is likely that
both MrGARK and MRI-GARK will be easier to extend to fifth-order and
to incorporate implicitness at the slow time scale.

We note that there are many avenues for further extensions of this
work.  On a theoretical level, we have yet to explore the 5th-order
conditions for RMIS methods, an effort that will undoubtedly leverage
Lemma \ref{lem:RMIS_bf_simplifications} repeatedly, and likely result
in a set of additional conditions on the outer table $T_O$.
Similarly, we note that although all of the theory presented in
Sections \ref{sec:MISGARK} and \ref{ch:RMIS} assume an autonomous ODE,
the inverter-chain problem is non-autonomous and still shows the
predicted orders of accuracy for all methods tested, hinting that this
theory can be extended to the non-autonomous context as well.  On a
purely computational level, we plan to investigate mixed
implicit/explicit multirate methods based on the RMIS structure
(implicit $T_I$ with EDIRK or ESDIRK structure), including selection
of $T_O$ and $T_I$ pairs for optimal efficiency.  There are also
numerous algorithmic extensions of this work, including (a)
exploration of the RMIS/MIS embedding to perform adaptivity (in both
$h$ and $m$) for efficient, tolerance-based, calculations; (b)
utilization of the telescopic property for three-rate problems with
large time-scale separations, and for $n$-rate problems with smaller
time-scale separations (as arise in explicit methods for hyperbolic
PDEs posed on spatially adaptive meshes); and (c) exploration of
RMIS-based ``self-adjusting'' methods \cite{Fok2016, Savcenco2007,
  Savcenco2008, Hundsdorfer2009, Savcenco2010}, that use the temporal
error estimate to automatically determine a fast/slow partitioning for
the problem.

\section*{Acknowledgements}
Support for this work was provided by the Department of Energy, Office
of Science project ``Frameworks, Algorithms and Scalable Technologies
for Mathematics (FASTMath),'' under Lawrence Livermore National
Laboratory subcontracts B598130, B621355, and B626484.

\appendix
\section{Extended Proofs}
\label{sec:Proofs}

\setcounter{theorem}{0}
\setcounter{lemma}{0}

\begin{theorem}
\notag
Assume that the inner base method $T_I$ is at least third order, the
outer base method $T_O$ is at least fourth order, and that both
satisfy the row-sum consistency conditions
\eqref{eq:consistency_conditions}.  If $T_O$ is explicit and satisfies
the additional condition \eqref{eq:rmis_4th}, i.e.,
\begin{equation*}
  v^{O\intercal} A^O c^O = \frac{1}{12},
\end{equation*}
where
\[
  v^O_i = \begin{cases}
    0,& i=1,\\
    b_i^O\(c_i^O - c_{i-1}^O\) + \(c_{i+1}^O - c_{i-1}^O\)\sum_{j=i+1}^{s^O}b_j^O,&1<i<s^O,\\
    b_{s^O}^O\(c_{s^O}^O-c_{s^O-1}^O\),&i=s^O,
  \end{cases}
\]
then the MIS and RMIS coefficients $\lam{A}{f,f}$, $\lam{A}{f,s}$,
$\lam{A}{s,f}$, $\lam{A}{s,s}$ and $\lam{b}{s}$ will satisfy all of
the ``slow'' fourth-order conditions
(i.e.,~\eqref{eq:mvOrderGARK4a}-\eqref{eq:mvOrderGARK4d} with
$\sigma=s$ and $\nu,\mu=\{f,s\}$).\\

Since $\lam{A}{s,s} = A^O$ and $\lam{b}{s} = b^O$, the equations
\eqref{eq:mvOrderGARK4a}-\eqref{eq:mvOrderGARK4d} with
$\sigma=\nu=\mu=s$ follow directly from assuming that $T_O$
is fourth-order.  The remaining fourth order conditions are
\begin{align*}
  \(\lam{b}{s}\times\lam{c}{s}\)^{\intercal}\lam{A}{s,f}\lam{c}{f} = \frac18,
    &\qquad
    \lamT{b}{s}\lam{A}{s,f}\(\lam{c}{f}\times\lam{c}{f}\) = \frac{1}{12},\\
  \lamT{b}{s}\vect{A}^{\{s,s\}}\lam{A}{s,f}\lam{c}{f} = \frac{1}{24},
    &\qquad
    \lamT{b}{s}\lam{A}{s,f}\vect{A}^{\{f,f\}}\lam{c}{f} = \frac{1}{24},\quad\text{and}\\
  \lamT{b}{s}\lam{A}{s,f}\vect{A}^{\{f,s\}}\lam{c}{s} = \frac{1}{24}&.
\end{align*}
We examine these in order:
\begin{align*}
  \(\lam{b}{s}\times\lam{c}{s}\)^{\intercal}\mathbf{A}^{\{s,f\}}\mathbf{c}^{\{f\}}
  &= \(\lam{b}{s}\times\lam{c}{s}\)^{\intercal}\(\frac12\lam{c}{s}\times\lam{c}{s}\)\\
  &= \frac12 \lamT{b}{s}\(\lam{c}{s}\times\lam{c}{s}\times\lam{c}{s}\)
  = \frac12 \(\frac14\) = \frac18,
\end{align*}
where we have used the identity $\mathbf{A}^{\{s,f\}}\mathbf{c}^{\{f\}}=\frac12\(\lam{c}{s}\)^2$
\cite[Theorem 3.1]{Sandu2013multirate}, and that $T_O$ satisfies \eqref{eq:mvOrderGARK4a}.
Similarly,
\small
\begin{align*}
  &\lamT{b}{s}\lam{A}{s,f}\(\lam{c}{f}\times\lam{c}{f}\)\\
  &= b^{O\intercal} \begin{bmatrix}
    c_2^O \vect{g}_2 b^{I\intercal} & \cdots &
    \(c_{s^O}^O-c_{s^O-1}^O\)\vect{g}_{s^O} b^{I\intercal} &
    \zeros \end{bmatrix}
  \begin{bmatrix} \(c_2^O c^I\)^2\\
    \(c_2^O\ones^{\{s^I\}}+\(c_3^O - c_2^O\) c^I\)^2\\
    \vdots\\
    \(c_{s^O}^O\ones^{\{s^I\}}+\(1-c_{s^O}^O\) c^I\)^2
  \end{bmatrix}\\
  &= b^{O\intercal} \bigg[
    \frac13 \(c_2^O\)^3 \vect{g}_2 + \(c_3^O-c_2^O\) \(\(c_2^O\)^2 + c_2^O\(c_3^O - c_2^O\) + \frac13 \(c_3^O - c_2^O\)^2 \)\vect{g}_3  + \cdots\\
    &\qquad\qquad +
    \(c_{s^O}^O-c_{s^O-1}^O\) \(\(c_{s^O-1}^O\)^2 + c_{s^O-1}^O\(c_{s^O}-c_{s^O-1}^O\) + \frac13 \(c_{s^O}-c_{s^O-1}^O\)^2 \)\vect{g}_{s^O} \bigg]\\
  &= \frac13 b^{O\intercal} \bigg[
     \(c_2^O\)^3 \vect{g}_2 + \( \(c_3^O\)^3 - \(c_2^O\)^3 \)\vect{g}_3 + \cdots +
     \( \(c_{s^O}^O\)^3 - \(c_{s^O-1}^O\)^3 \)\vect{g}_{s^O} \bigg]\\
  &= \frac13 b^{O\intercal}
    \begin{bmatrix} 0\\ \(c_2^O\)^3 \\ \vdots \\ \(c_{s^O}^O\)^3 \end{bmatrix}
  = \frac13 b^{O\intercal} \(c^O\times c^O\times c^O\) = \frac13 \(\frac14\)
    = \frac{1}{12},
\end{align*}
\normalsize
where we have relied on the fact that $T_I$ is third order, and that
$T_O$ is explicit and fourth order. Again using the result
$\mathbf{A}^{\{s,f\}}\mathbf{c}^{\{f\}}=\frac12\(\lam{c}{s}\)^2$,
along with the fact that $T_O$ is fourth-order, we have
\begin{align*}
  \lamT{b}{s} \vect{A}^{\{s,s\}} \lam{A}{s,f} \lam{c}{f}
  &= \lamT{b}{s} \vect{A}^{\{s,s\}} \(\frac12\lam{c}{s}\times \lam{c}{s}\)\\
  &= \frac12 \(b^O\)^\intercal A^O \(c^O\times c^O\)
  = \frac12 \(\frac{1}{12}\) = \frac{1}{24}.
\end{align*}
Similarly,
\footnotesize
\begin{align*}
  &\lamT{b}{s} \lam{A}{s,f} \vect{A}^{\{f,f\}} \lam{c}{f}\\
  &= b^{O\intercal} \lam{A}{s,f}
    \begin{bmatrix}
      c_2^O A^I & \zeros & \cdots & \zeros\\
      c_2^O \ones^{\{s^I\}} b^{I\intercal} & \(c_3^O - c_2^O\)A^I & \cdots & \zeros\\
      \vdots &  & \ddots & \\
      c_2^O\ones^{\{s^I\}} b^{I\intercal} & \(c_3^O - c_2^O\)\ones^{\{s^I\}} b^{I\intercal} & \cdots & \(1-c_{s^O}^O\)A^I
    \end{bmatrix} \begin{bmatrix}
      c_2^O c^I\\
      c_2^O\ones^{\{s^I\}}+\(c_3^O - c_2^O\) c^I\\
      \vdots\\
      c_{s^O}^O\ones^{\{s^I\}}+\(1-c_{s^O}^O\) c^I
    \end{bmatrix}\\
  &= b^{O\intercal} \begin{bmatrix}
      c_2^O \vect{g}_2 b^{I\intercal} & \cdots &
      \(c_{s^O}^O-c_{s^O-1}^O\)\vect{g}_{s^O} b^{I\intercal} &
      \zeros
    \end{bmatrix} \begin{bmatrix}
      \(c_2^O\)^2 A^I c^I\\
      \frac12 \(c_2^O\)^2 \ones^{\{s^I\}} + \(c_3^O - c_2^O\)\(c_2^Oc^I+\(c_3^O - c_2^O\) A^I c^I\)\\
      \vdots\\
      \frac12 \(c_{s^O}^O\)^2 \ones^{\{s^I\}} + \(1-c_{s^O}^O\)\(c_{s^O}^Oc^I+\(1-c_{s^O}^O\) A^I c^I\)
    \end{bmatrix}\\
  &= b^{O\intercal} \sum_{i=2}^{s^O} \vect{g}_i\left[ \(c_{i}^O-c_{i-1}^O\) b^{I\intercal} \(\frac12
    \(c_{i-1}^O\)^2 \ones^{\{s^I\}} + \(c_i^O - c_{i-1}^O\)\(c_{i-1}^Oc^I+\(c_i^O - c_{i-1}^O\) A^I c^I\)\) \right]\\
  &= b^{O\intercal} \sum_{i=2}^{s^O} \vect{g}_i \left[ \(c_{i}^O-c_{i-1}^O\) \(\frac12 \(c_{i-1}^O\)^2b^{I\intercal} \ones^{\{s^I\}} + \(c_i^O - c_{i-1}^O\)\(c_{i-1}^Ob^{I\intercal} c^I+\(c_i^O - c_{i-1}^O\)b^{I\intercal} A^I c^I\)\) \right]\\
  &= b^{O\intercal} \sum_{i=2}^{s^O} \vect{g}_i \left[ \(c_{i}^O-c_{i-1}^O\) \(\frac12
    \(c_{i-1}^O\)^2 + \(c_i^O -  c_{i-1}^O\)\(\frac12 c_{i-1}^O + \frac16 \(c_i^O - c_{i-1}^O\) \)\) \right]\\
  &= \frac16 b^{O\intercal} \sum_{i=2}^{s^O} \vect{g}_i \left[ \(c_i^O\)^3 - \(c_{i-1}^O\)^3 \right]\\
  &= \frac16 b^{O\intercal} \begin{bmatrix}
      0\\
      \(\(c_2^O\)^3 - \(c_1^O\)^3\)\\
      \(\(c_3^O\)^3 - \(c_2^O\)^3\) + \(\(c_2^O\)^3 - \(c_1^O\)^3\)\\
      \vdots\\
      \(\(c_{s^O}^O\)^3 - \(c_{s^O-1}^O\)^3\) + \cdots + \(\(c_2^O\)^3 - \(c_1^O\)^3\)
    \end{bmatrix}
  = \frac16 b^{O\intercal} \begin{bmatrix}
      \(c_1^O\)^3\\
      \(c_2^O\)^3\\
      \(c_3^O\)^3\\
      \vdots\\
      \(c_{s^O}^O\)^3
    \end{bmatrix}\\
  &= \frac16 b^{O\intercal} \(c^O\times c^O\times c^O\)
  = \frac16 \(\frac14\) = \frac{1}{24},
\end{align*}
\normalsize
where we have used the fact that $T_I$ is third order, and $T_O$ is both
fourth order and explicit.\\

Our final slow fourth-order condition becomes
\footnotesize
\begin{align*}
  &\lamT{b}{s} \lam{A}{s,f} \lam{A}{f,s} \lam{c}{s}\\
  &= b^{O\intercal} \begin{bmatrix}
    c_2^O \vect{g}_2 b^{I\intercal} & \cdots &
    \(c_{s^O}^O-c_{s^O-1}^O\)\vect{g}_{s^O} b^{I\intercal} & \zeros
  \end{bmatrix} \begin{bmatrix}
    c^I \vect{e}_2^{\intercal} A^O\\
    \vdots\\
    \ones^{\{s^I\}} \vect{e}_{i}^{\intercal} A^O + c^I \(\vect{e}_{i+1}-\vect{e}_{i}\)^{\intercal} A^O\\
    \vdots\\
    \ones^{\{s^I\}} \vect{e}_{s^O}^{\intercal} A^O + \vect{c}^I \(\vect{b}^{O\intercal}-\vect{e}_{s^O}^{\intercal} A^O\)
  \end{bmatrix} c^O \\
  &= \frac12 b^{O\intercal}\(c_2^O \vect{g}_2 \vect{e}_2^{\intercal} + \sum_{i=2}^{s^O-1}  \(c_{i+1}^O-c_{i}^O\)\vect{g}_{i+1} \( \vect{e}_{i+1} + \vect{e}_{i} \)^{\intercal}\) A^O c^O \\
  &= \frac12 b^{O\intercal}\begin{bmatrix}
    \zeros^{\intercal}\\
    c_2^O\vect{e}_2^{\intercal} \\
    c_3^O\vect{e}_2^{\intercal} + \(c_3^O - c_2^O\)\vect{e}_3^{\intercal}\\
    c_3^O\vect{e}_2^{\intercal} + \(c_4^O - c_2^O\)\vect{e}_3^{\intercal} + \(c_4^O-c_3^O\)\vect{e}_4^{\intercal}\\
    \vdots\\
    c_3^O\vect{e}_2^{\intercal} + \(c_4^O - c_2^O\)\vect{e}_3^{\intercal} + \cdots + \(c_{s^O-1}^O-c_{s^O-2}^O\)\vect{e}_{s^O-1}^{\intercal}\\
    c_3^O\vect{e}_2^{\intercal} + \(c_4^O - c_2^O\)\vect{e}_3^{\intercal} + \cdots + \(c_{s^O}^O-c_{s^O-2}^O\)\vect{e}_{s^O-1}^{\intercal} + \(c_{s^O}^O-c_{s^O-1}^O\)\vect{e}_{s^O}^{\intercal}
  \end{bmatrix} A^O c^O \\
  &= \frac12 \begin{bmatrix}
    0\\
    b_2^Oc_2^O + \(c_3^O-c_1^O\)\sum_{i=3}^{s^O}b_i^O \\
    b_3^O\(c_3^O - c_2^O\) + \(c_4^O - c_2^O\)\sum_{i=4}^{s^O}b_i^O \\
    \vdots\\
    b_{s^O-1}^O\(c_{s^O-1}^O-c_{s^O-2}^O\) + b_{s^O}^O\(c_{s^O}^O-c_{s^O-2}^O\) \\
    b_{s^O}^O\(c_{s^O}^O-c_{s^O-1}^O\)
  \end{bmatrix}^{\intercal} A^O c^O
  = \frac12 v^{O\intercal} A^O c^O = \frac{1}{24},
\end{align*}
\normalsize
where we have relied on the assumption \eqref{eq:rmis_4th} and the
fact that $T_I$ is at least second-order accurate.
\qed
\end{theorem}

\begin{lemma}
\notag
Suppose that the coefficients $\lam{b}{f}$ are chosen as in equation
\eqref{eq:RMIS_bf}, and that the inner Butcher table $T_I$ has
explicit first stage (i.e.~the first entry of $c^I$ and the first row
of $A^I$ are identically zero).  Then the following identities hold:
\begin{align}
  \notag
  \lamT{b}{f}\(\lam{c}{f}\)^q &= \lamT{b}{s} \(\lam{c}{s}\)^q,
                                \;\forall q\ge 0,\\
  \notag
  \lamT{b}{f} \lam{A}{f,f} &= \lamT{b}{s} \lam{A}{s,f},\\
  \notag
  \lamT{b}{f} \lam{A}{f,s} &= \lamT{b}{s} \lam{A}{s,s},\\
  \(\lam{b}{f}\times\lam{c}{f}\)^{\intercal} \lam{A}{f,f} &= \(\lam{b}{s}\times\lam{c}{s}\)^{\intercal} \lam{A}{s,f},\\
  \(\lam{b}{f}\times\lam{c}{f}\)^{\intercal} \lam{A}{f,s} &= \(\lam{b}{s}\times\lam{c}{s}\)^{\intercal} \lam{A}{s,s},
\end{align}
where $\lam{A}{f,f}$, $\lam{A}{s,f}$ and $\lam{A}{f,s}$ are defined as
in equations \eqref{eq:MIS_Aff}-\eqref{eq:MIS_Afs}, $\lam{c}{f}$ is
defined as in \eqref{eq:MIS_cf}, $\lam{c}{s} = c^O$, and $\lam{A}{s,s} = A^O$.\\

These follow from direct application of the vector-matrix products:
\begin{align*}
  &\lamT{b}{f}\(\lam{c}{f}\)^q =
  \begin{bmatrix}
    b_1^O \vect{e}_1^{\intercal} & b_2^O \vect{e}_1^{\intercal} &
    \cdots & b_{s^O}^O \vect{e}_1^{\intercal}
  \end{bmatrix} \begin{bmatrix}
    \(c_2^O c^I\)^q\\
    \(c_2^O\ones^{\{s^I\}}+\(c_3^O - c_2^O\) c^I\)^q\\
    \vdots\\
    \(c_{s^O}^O\ones^{\{s^I\}}+\(1-c_{s^O}^O\) c^I\)^q
  \end{bmatrix} \\
  &= b^{O\intercal} \(c^O\)^q = \lamT{b}{s} \(\lam{c}{s}\)^q,\\
  &\lamT{b}{f} \lam{A}{f,f} \\
  &= \lamT{b}{f} \begin{bmatrix}
    c_2^O A^I & \zeros & \cdots & \zeros\\
    c_2^O \ones^{\{s^I\}} b^{I\intercal} & \(c_3^O - c_2^O\)A^I & \cdots & \zeros\\
    \vdots &  & \ddots & \\
    c_2^O\ones^{\{s^I\}} b^{I\intercal} & \(c_3^O - c_2^O\)\ones^{\{s^I\}} b^{I\intercal} & \cdots & \(1-c_{s^O}^O\)A^I
  \end{bmatrix}\\
  &= \begin{bmatrix}
    \( b_2^O + \cdots + b_{s^O}^O \) c_2^O b^{I\intercal} &
    \( b_3^O + \cdots + b_{s^O}^O \) \(c_3^O - c_2^O\) b^{I\intercal} &
    & \cdots & \zeros^{\intercal}
  \end{bmatrix}\\
  &= b^{O\intercal}
  \begin{bmatrix}
    c_2^O \vect{g}_2 b^{I\intercal} & \cdots &
    \(c_{s^O}^O-c_{s^O-1}^O\)\vect{g}_{s^O} b^{I\intercal} &
    \zeros
  \end{bmatrix} = \lamT{b}{s} \lam{A}{s,f},
\end{align*}
and
\begin{align*}
  &\lamT{b}{f} \lam{A}{f,s} \\
  &= \begin{bmatrix}
    b_1^O \vect{e}_1^{\intercal} & b_2^O \vect{e}_1^{\intercal} &
    \cdots & b_{s^O}^O \vect{e}_1^{\intercal}
  \end{bmatrix} \begin{bmatrix}
    c^I \vect{e}_2^{\intercal} A^O\\
    \vdots\\
    \ones^{\{s^I\}}\vect{e}_{i}^{\intercal} A^O + c^I \(\vect{e}_{i+1}-\vect{e}_{i}\)^{\intercal} A^O\\
    \vdots\\
    \ones^{\{s^I\}}\vect{e}_{s^O}^{\intercal} A^O + c^I \(\vect{b}^{O\intercal}-\vect{e}_{s^O}^{\intercal} A^O\)
  \end{bmatrix}\\
  &= \sum_{i=2}^{s^O} b_i^O \vect{e}_{i}^{\intercal} A^O =
  b^{O\intercal} A^O = \lamT{b}{s} \lam{A}{s,s}.
\end{align*}
Proof of the final two conditions are essentially identical to those
above, by using the substitutions $\lam{b}{f} \to \(\lam{b}{f}\times\lam{c}{f}\)$ and
$\lam{b}{s} \to \(\lam{b}{s}\times\lam{c}{s}\)$.
\qed
\end{lemma}

\section{Butcher table information}
\label{sec:Coefficients}
Butcher's general solution for a 4-stage fourth order explicit Runge-Kutta
method depends only on two free variables, $c_2$ and $c_3$, and is
given by \cite{Butcher2008}
\begin{align}
  a_{2,1} &= c_2, \notag\\
  a_{3,1} &= \frac{c_3\(c_3+4c_2^2-3c_2\)}{2c_2\(2c_2-1\)},
    \notag\\
  a_{3,2} &= -\frac{c_3\(c_3-c_2\)}{2c_2\(2c_2-1\)},
    \notag\\
  a_{4,1} &= \frac{-12c_3c_2^2+12c_3^2c_2^2+4c_2^2-6c_2+15c_2c_3-12c_3^2c_2+2+4c_3^2-5c_3}{2c_2c_3\(-4c_3+6c_3c_2+3-4c_2\)},
  \notag\\
  a_{4,2} &= \frac{\(c_2-1\)\(4c_3^2-5c_3+2-c_2\)}{2c_2\(c_3-c_2\)\(-4c_3+6c_3c_2+3-4c_2\)}, \notag\\
  a_{4,3} &= -\frac{\(2c_2-1\)\(c_2-1\)\(c_3-1\)}{c_3\(c_3-c_2\)\(-4c_3+6c_3c_2+3-4c_2\)}, \label{eq:RK4stage}\\
  b_1 &= \frac{6c_3c_2-2c_3-2c_2+1}{12c_3c_2}, \notag\\
  b_2 &= -\frac{\(2c_3-1\)}{12c_2\(c_2-1\)\(c_3-c_2\)}, \notag\\
  b_3 &= \frac{\(2c_2-1\)}{12c_3\(c_2-c_3c_2+c_3^2-c_3\)}, \notag\\
  b_4 &= \frac{-4c_3+6c_3c_2+3-4c_2}{12\(c_3-1\)\(c_2-1\)}. \notag
\end{align}
For $T_O$ with this structure, the 3rd-order MIS condition
\eqref{eq:rfsmr_3rd} is equivalent to
\begin{align}
\label{eq:RK4stageSolRFSMR}
  3(c_2-1) (6c_2^{2}c_3^{2}-4c_2^{2}c_3-6c_2c_3^{3}&+8c_2c_3^{2}-11c_2c_3+6c_2+4c_3^{3}-7c_3^{2}+7c_3-3) \notag\\
  & - 2(2c_2-1)(4c_2+4c_3-6c_2c_3-3)=0,
\end{align}
and the 4th-order RMIS condition \eqref{eq:rmis_4th} is equivalent to
\begin{align}
  \label{eq:RK4stageSolRMIS}
  36c_3^{4}-120c_3^{3}+80c_3^{2}-12c_3+1 - \(4c_2(3c_3+1) - 6c_3^{2} + 2c_3 - 3\)^2 = 0.
\end{align}



\bibliography{rmis_paper}

\end{document}